\documentclass{article}
\usepackage{arxiv}
\usepackage{graphicx,amsmath,amsfonts,amscd,amsthm,amssymb,pstricks}
\usepackage{tabularx}
\newcolumntype{Y}{>{\centering\arraybackslash}X}
\usepackage{booktabs}
\usepackage{color}% Include figure files
\usepackage{dcolumn}% Align table columns on decimal point
\usepackage{bm}% bold math
\usepackage{longtable}
\usepackage{mathrsfs}
\usepackage{textcomp}
\usepackage{esvect}
\usepackage{float}
\usepackage[USenglish,british]{babel}
\usepackage{environ}
\usepackage{xifthen}
\usepackage{xargs}		% programming: better newcommand
\usepackage{hyperref}
\usepackage[separate-uncertainty = true]{siunitx}
\usepackage{physics,mathtools}
\usepackage{multirow}
\usepackage{comment}
\usepackage{enumerate}
\usepackage{appendix}
\usepackage{tikz}	
\usetikzlibrary{arrows, calc, matrix, patterns, decorations.markings, shapes, decorations.pathmorphing, shadows.blur}
\usepackage[utf8]{inputenc}
\usepackage{enumitem}
\usepackage{bm}% bold math
\def\Tr{\hbox{Tr}}
\def\det{\hbox{det}}

\font\grassettogreco=cmmib10
\font\scriptgrassettogreco=cmmib7
\font\scriptscriptgrassettogreco=cmmib10 at 5 truept
\font\sansserif=cmss10
\font\scriptsansserif=cmss10 at 7 truept
\font\scriptscriptsansserif=cmss10 at 5 truept
\def\ssm{\fam=14}
\input amssym.def
\input amssym.tex
\def \ebf{{\bf e}}

\def \xbf{{\bf x}}
\def \ybf{{\bf y}}

\def \Iop{{\mathchardef\alpha="710B \ssm \char'111}}

\def \Lop{{\mathchardef\alpha="710B \ssm \char'114}}

\def \Rop{{\mathchardef\alpha="710B \ssm \char'122}}

\def \Wop{{\mathchardef\alpha="710B \ssm \char'127}}
\def \Xop{{\mathchardef\alpha="710B \ssm \char'130}}

\def \eps{\epsilon}

  2
  3
 4
 5
   \newcount\Silicon \global\Silicon=1
   \newcount\PC \global\PC=2
   \newcount\parziale
   \def\Colori#1{\global\parziale=#1
                \ifnum\parziale=\Silicon
                    \input colors
                    \gdef\Color##1{\Black{##1}}                    
                \else\ifnum\parziale=\PC
                    \input colordvi
                    \gdef\textRGB##1{\textColor{##1 0.}}
                    \gdef\GrayA##1{\textGray{##1}}
                    \gdef\GrayB##1{\textGray{##1}}
                    \gdef\GrayC##1{\textGray{##1}}
                    \gdef\GrayD##1{\textGray{##1}}
                    \gdef\GrayE##1{\textGray{##1}}
                    \gdef\GrayF##1{\textGray{##1}}
                    \gdef\GrayG##1{\textGray{##1}}
                    \gdef\GrayH##1{\textGray{##1}}\fi                   
                \fi}        
\def\mean#1{\langle \,#1\,\rangle}

\def \parton#1{\left({{#1}}\right)}

\def \parmean#1{\left\langle{{#1}}\right\rangle}
\def \underwrite#1{\mathop{\vtop{\ialign {##\crcr
$\hfil\displaystyle {#1}\hfil$\crcr\noalign{\kern3pt\nointerlineskip}
\crcr\noalign{\kern3pt}}}} \limits}
\def\sqr#1#2{{\vcenter{\hrule height.#2pt
     \hbox{\vrule width.#2pt height#1pt \hskip#1pt
       \vrule width.#2pt}
     \hrule height.#2pt}}}

\def\DAL{\hbox{\raise.250ex \hbox{$\sqr7{10}\,$}}} % DALambertiano

\def\bib#1.#2/#3/#4/#5/#6/#7.{\frenchspacing\item{[#1]}#2:\ {\it ``#3''}
~--~#4\ $\underline{\bf #5}$,\ #6 (#7)}
\def\diagramma#1#2#3#4#5#6#7#8{
 \vbox to 2.5cm{
       \hbox to 3cm{\hfil ${#1}$ \hfil}
       \hbox to 3cm{ ${#2}$\rightarrowfill ${#3}$ }
\hbox to 3cm{${#4} \Biggl\uparrow\hfill\Biggr\uparrow {#5}$}
       \hbox to 3cm{ ${#6}$\rightarrowfill ${#7}$ }
       \hbox to 3cm{\hfil ${#8}$ \hfil}   }    }
\def\sopra#1#2{{\raise 0.8 ex
\hbox{$
{{\scriptstyle \,{#2}}	\atop \displaystyle{#1}}
$}}
}
\def\figuraps#1#2#3#4#5{
\par
\midinsert
\centerline{\bf #4}
\vbox to #3 truecm{
\vskip #3 truecm
\ifnum #1 = 0	% ---  siamo su PC
\special {ps: plotfile #2}
\else		% ---  siamo su MVX
\special {#2 0 0 moveto 16} \fi
}
\centerline{#5}
\endinsert}
\font\cofon = cmr6
\font\cobfon = cmbx6
\font\copi = cmr9
\def\codlib{{\copi\copyright}{\cofon 88-08- }{\cobfon 9820}}
\def\riga{\vskip .1  truecm   \hrule \vskip .2	    truecm \noindent }
\def\HeadLinea#1#2{	\headline={\vbox to 0pt{\vss\noindent
{\ifnum \pageno=1  \hfill {\bf \folio}		 %  prima pagina
\else {\ifodd \pageno			   % pagina sinistra
{\noindent     \hfill  {\it     #2} \quad {\bf \folio}
 }\riga
\else				 %pagina destra
{\noindent {\bf\folio} \quad  {\it  #1} \hfill 	}
\riga
\fi } \fi }			}}
}
\def\testatina#1#2{	\headline={\vbox to 0pt{\vss\noindent
{\ifnum \pageno=1  \hfill {\bf \folio}		 %  prima pagina
\else {\ifodd \pageno			   % pagina sinistra
{\noindent  \codlib   \hfill  {\it     #2} \quad {\bf \folio}
 }\riga
\else				 %pagina destra
{\noindent {\bf\folio} \quad  {\it  #1} \hfill \codlib	}
\riga
\fi } \fi }			}}
}
\def\testatinacap#1#2#3{	\headline={\vbox to 0pt{\vss\noindent
{\ifnum \pageno=#3  \hfill {\bf \folio}		 %  prima pagina
\else {\ifodd \pageno			   % pagina sinistra
{\noindent  \codlib   \hfill  {\it     #2} \quad {\bf \folio}
 }\riga
\else				 %pagina destra
{\noindent {\bf\folio} \quad  {\it  #1} \hfill \codlib	}
\riga
\fi } \fi }			}}
}
\def\oggi{\number\day\space\ifcase\month
   \or gennaio\or febbraio\or marzo\or aprile\or maggio\or giugno\or
   luglio\or agosto\or settembre\or ottobre\or novembre\or dicembre
   \fi\space\number\year}
\def\today{\number\day\space\ifcase\month
   \or January\or February\or March\or April\or May\or June\or
   July\or August\or September\or October\or November \or December
   \fi\space\number\year}   
\def\frame#1{\ifmmode\dframe{#1}\else\leavevmode\lower 2.4 pt
    \hbox{\vrule\unskip\vbox{\hrule\kern 1.5 pt\hbox{\kern
    1.5 pt{#1}\kern 0.5 pt}\kern 2 pt\hrule}\unskip\vrule}\fi}

\def\dframe#1{\hbox{\vrule\unskip$\vcenter{\hrule\kern 3 pt\hbox
    {\kern 3 pt$\displaystyle{#1}$\kern3pt}\kern 3 pt\hrule}$\vrule}}

\title{Performance analysis of indicators of chaos for nonlinear dynamical systems}

\author{A. Bazzani\\
Dipartimento di Fisica e Astronomia, Universit\`a di Bologna and INFN Bologna, via Irnerio 46, Bologna, Italy\\
\And
M. Giovannozzi\thanks{Corresponding author: massimo.giovannozzi@cern.ch}\\
Beams Department, CERN, Esplanade\ des Particules 1, 1211 Geneva 23, Switzerland\\
\And
C.E. Montanari\\
Dipartimento di Fisica e Astronomia, Universit\`a di Bologna and INFN Bologna, via Irnerio 46, Bologna, Italy\\
Beams Department, CERN, Esplanade\ des Particules 1, 1211 Geneva 23, Switzerland\\
\And
G. Turchetti\\
Dipartimento di Fisica e Astronomia, Universit\`a di Bologna and INFN Bologna, via Irnerio 46, Bologna, Italy}

\begin{document}
\maketitle

\begin{abstract}
The efficient detection of chaotic behavior in orbits of a complex dynamical system is an active domain of research. Several indicators have been proposed in the past, and new ones have recently been developed in view of improving the performance of chaos detection by means of numerical simulations. The challenge is to predict chaotic behavior based on the analysis of orbits of limited length. In this paper, the performance analysis of past and recent indicators of chaos, in terms of predictive power, is carried out in detail using the dynamical system characterized by a symplectic H\'enon-like cubic polynomial map.  
\end{abstract}

\section{\label{sec:intro}Introduction}
The study of the long-term evolution of Hamiltonian systems is a very difficult task from both a theoretical and a numerical point of view. The KAM theory\cite{KAM4} does not provide a solution to the stability problem for Hamiltonian systems in more than two degrees of freedom. Therefore, great effort has been devoted to improving time stability estimates after the celebrated Nekhoroshev theorem~\cite{Nekhoroshev:1977aa}. However, the existence of chaotic layers in phase space strongly affects the long-term evolution of the orbits, and for this reason, numerical indicators have been proposed to detect the chaotic character of orbits using a limited number of time steps. 

For a given Hamiltonian model, one has to tackle the problem of comparing the performance of the various indicators to assess which one provides the optimal classification of the orbits. In applications, this task must be accomplished taking into account the characteristics of the physical problem under consideration. For instance, in the field of accelerator physics, the study of the charged-hadron motion in the magnetic lattice of a circular accelerator is often devoted to the determination of the region of phase space in which bounded motion occurs. The extent of such a region is called dynamic aperture, and its precise determination involves studying the stability of orbits of a 6D symplectic map in a neighborhood of an elliptic fixed point, up to $10^8-10^9$ iterations~(see, e.g., \cite{Bazzani:262179}). An exhaustive analysis of the phase-space topology is clearly beyond the current computational capabilities, even for relatively simple systems. Therefore, indicators of chaos turn out to be extremely useful to reduce the amount of computational time needed to assess the character of orbits (regular or chaotic). This task may be affected by the presence of orbit diffusion in phase space, which occurs in chaotic layers. The presence of small stochastic effects, which naturally arise in physical systems, may prevent orbit trapping near regular regions, the so-called stickiness phenomenon~\cite{Kandrup1999DiffusionAS,SZEZECH2005394}, thus inducing diffusive behavior in phase space. 

It is worth noting that polynomial symplectic maps are central for the analysis of accelerator physics problems, but they are also present in other domains and have been intensively studied to understand the phase space structure of Hamiltonian systems~\cite{Turaev2003}, and are a fundamental tool for long-term integration of orbits~\cite{Bazzani:262179}. 

The main result of this paper is to show that it is possible to determine a classification performance ranking of the main commonly used chaotic indicators when applied to a generic $4d$ cubic polynomial symplectic map of H\'enon-like form (see, e.g., ~\cite{Bazzani:262179}), which is an excellent prototype dynamical system for applications, such as circular hadron accelerators.

The indicators of chaos are typically based on the existence of positive Lyapunov characteristic exponents, and their numerical performance is strongly affected in the regions where sticky orbits are present.

The family of Fast Lyapunov indicators ($FLI$)~\cite{Froeschle1997} has been proposed to distinguish the regions of regular and chaotic motion for symplectic maps~\cite{SZEZECH2005394}. They also proved to be suitable for identifying resonant regions in phase space and to visualize the Arnold web of resonances where slow diffusion occurs~\cite{Arnold:937549}. These indicators are based on the evolution of an initial deviation vector and provide the linear response of the tangent map along an orbit. When considering one or more initial deviation vectors, the result depends on the direction of the initial deviation vectors. To overcome this, the linear response to a random displacement vector with zero mean value and unit variance was recently proposed~\cite{Turchetti2017}. The trace of the corresponding covariance matrix defines the square Lyapunov Error ($LE$), which is similar to $FLI$. Furthermore, the invariants of the covariance matrix of order $k>1$ are asymptotically related to the sum of the first $k$ Lyapunov exponents. However, unlike the Generalized Alignment Index ($GALI^{(k)}$) indicators~\cite{Bountis2007,Skokos2015}, these invariants do not depend on the initial deviations~\cite{Skokos2008}. Recently, a couple of approaches have been proposed to improve the performance of some indicators, namely applying the Weighted Birkhoff averaging~\cite{Das_2017} or the Mean Exponential Growth of Nearby Orbit ($MEGNO$)~\cite{Cincotta1996}, which is used to filter the oscillations and to improve the accuracy by averaging on map iterations~\cite{Gozdziewski01,Mestre2011}.

To calculate the sensitivity to small deviations along an orbit, the Reversibility Error Method ($REM$) can be used~\cite{Panichi2016,Panichi2017}. In this case, the linear response to the forward evolution in the presence of small random noise is considered, followed by the unperturbed backward evolution. The covariance matrix of the random process, which provides the final deviation from the initial condition in the limit of zero noise amplitude, can be computed, and its invariants quantify the violation of reversibility. The first invariant for the forward-backward process $BF$ is the square of the reversibility error, which is equal to the sum of the squares of Lyapunov errors computed at each iteration of the map. This invariant can be compared with the results for $REM$, when the stochastic perturbation is generated by the finite numerical precision present in both the forward and backward directions. 

Finally, a completely different indicator introduced by J.Laskar~\cite{Laskar1999,Laskar2003} is represented by the Frequency Map Analysis $FMA$, which computes the variation of the main frequency of a given orbit considering different orbit lengths to detect the chaotic character. 

In this paper, we perform an accurate analysis of the performance of the indicators briefly introduced above to classify the orbits of a $4d$ modulated polynomial symplectic map, namely a $4d$ H\'enon map that is considered a reference model for several applications. In Section~\ref{sec:review} we define mathematically and discuss in some detail the chaos indicators considered, and in Section~\ref{sec:numerical_implementations} we discuss their numerical implementation. In Section~\ref{sec:results} we present the numerical results and rank the different indicators in terms of classification efficiency, in particular, studying their predictivity. Finally, some conclusions are drawn in Section~\ref{sec:conc}. In addition, we report some details on the computational cost of implementing indicators using parallel computing facilities in Appendix~\ref{app:computing}, while some considerations on the time dependence of indicators are presented in Appendix~\ref{app:timedep}. 
\section{Definition and main properties of indicators of chaos} \label{sec:review}
\subsection{Frequency Map Analysis\label{subsec:fma}}
Originally introduced by J.~Laskar in the field of celestial mechanics, the Frequency Map Analysis ($FMA$) rapidly found applications outside the initial domain of application (see, e.g., \cite{laskar1995frequency,lega1996numerical,papaphilippou1996frequency,papaphilippou1998global,Laskar1999, Papaphilippou1999, laskar2000application,PhysRevSTAB.4.124201,1288929,Papaphilippou:PAC03-RPPG007,Laskar2003,PhysRevSTAB.6.114801,shun2009nonlinear,PhysRevSTAB.14.014001,papaphilippou2014,tydecks:ipac18-mopmf057,PhysRevAccelBeams.22.071002} for a selected list of references, with special emphasis on accelerator-related applications) is a numerical method to inspect the global dynamics of multidimensional Hamiltonian systems, taking advantage of the quasiperiodicity of regular orbits located on KAM tori.

Given a Hamiltonian system $H(I,\theta) = H_0 (I) + \varepsilon H_1 (I, \theta)$, where for $\varepsilon=0$ the Hamiltonian $H_0(I)$ is integrable and $(I, \theta$ are action angle variables in $\mathbb{R}^n \times \mathbb{T}^n$, where $\mathbb{T}$ represents a one-dimensional torus. If the system is nondegenerate,
\begin{equation}
    \operatorname{det}\left(\frac{\partial \nu(I)}{\partial I}\right)=\operatorname{det}\left(\frac{\partial^2 H_0(I)}{\partial I^2}\right) \neq 0 \,,
\end{equation}
the application
\begin{equation}
    \begin{array}{r}
    F: I\in \mathbb{R}^{n} \longrightarrow \nu\in \mathbb{R}^n 
    \end{array}
\end{equation}
is a diffeomorphism on its image. This means that the invariant tori are equally identified by the action variables $I$ or by their corresponding frequency vector $\mathbf{\nu}$. For a nondegenerate system, when $\varepsilon$ is sufficiently small, the KAM theorem~\cite{KAM1,KAM2,KAM3}, states that there still exists a set of initial conditions of positive measure that correspond to regular orbits on invariant tori, for which, according to Pöschel~\cite{Poschel1982}, a similar diffeomorphism still applies.

Based on this theoretical framework, it is possible to distinguish between regular orbits on the KAM tori, which feature a discrete structure for Fourier components defined by the harmonic of the fundamental frequencies, and chaotic orbits, which exhibit a complex structure in the Fourier spectrum~\cite{1288929}. In this sense, $FMA$ is a technique that performs numerical evaluations of the frequency vector $\vb{\nu}$ from a time series corresponding to a certain interval $[i, i+n]$, for different values of $i$. In case of a regular orbit lying on a KAM tori, the frequency vectors for various $i$ will agree up to the precision of the numerical method used to determine the frequency. On the other hand, a chaotic orbit will have $\vb{\nu}$ that evolves over different intervals, showing fluctuations in frequency space~\cite{laskar2000application}.

To achieve an accurate numerical evaluation of fundamental frequencies, multiple studies have been carried out to improve standard algorithms such as the Fast Fourier Transform (FFT) or the Average Phase Advance (APA)~\cite{laskar1992measure,Laskar1999, Bartolini:316949,bartolini1998computer}. In the work of Bartolini et al.~\cite{Bartolini:292773}, the fundamental frequency is evaluated using an FFT combined with a Hanning filter and an interpolation algorithm, resulting in a closed-form formula for the fundamental frequency. In recent studies~\cite{russo:ipac2021-thpab189}, the frequency determination carried out using the average phase advance algorithm is improved by applying the weighted Birkhoff averaging~\cite{Das_2018}, which will be used in the sequel to perform the evaluation of $FMA$.
More precisely, we define $FMA_n$ as the Euclidean distance between two vectors defined by the fundamental frequencies $\nu_1$ and $\nu_2$, evaluated respectively over the time intervals $[0, n/2]$ and $[n/2, n]$ of the orbit. An initial condition on a KAM torus has $FMA_n$ converge to zero when $n\to \infty$. In contrast, an initial condition in a chaotic layer will converge to an asymptotic value for $FMA_n$ bounded away from zero.
\subsection{Lyapunov Error invariants\label{subsec:le}}
Let $M(\xbf,n)$ be a time-dependent symplectic map with $\xbf\in \mathbb{R}^{2d}$ where the first $d$ components of $\xbf$ are the space coordinates and the last $d$  their conjugate moments. Denoting by $DM$ the Jacobian matrix $(DM)_{ij}=\partial M_i/\partial x_j$ and by $\xbf_n$ the orbit after $n$ iterations, the corresponding tangent map $\Lop_n(\xbf)$ is defined by 
\begin{equation}
  \begin{aligned}
  &\xbf_n =M(\xbf_{n-1},n-1) \equiv M_n(\xbf) \,,   &\xbf_0=\xbf;& \\ \\
  &\Lop_n(\xbf) = DM(\xbf_{n-1}, n-1) \,\Lop_{n-1}(\xbf)  \equiv DM_n(\xbf) \,,     &\Lop_0=\Iop \, , &
  \end{aligned}
  \label{eq2-1}
\end{equation}
where $M_n(\xbf)= M(\xbf,n-1)\circ M_{n-1}(\xbf)$ with $M_0(\xbf)=\xbf$.

For any initial condition $\xbf$, consider a small stochastic deviation $\eps\boldsymbol{\xi}$ where
$\boldsymbol{\xi}$ is a unit random vector with $\mean{\boldsymbol{\xi}}=0$ and a unit covariance matrix $\mean{\boldsymbol{\xi}\,\boldsymbol{\xi}^T}=\Iop $, where the suffix $^T$ denotes the transposed vector. Letting $\ybf_n=M(\ybf_{n-1},n-1)$ be the orbit with initial condition $\ybf_0=\xbf_0+\eps\boldsymbol{\xi}$ the linear response $\boldsymbol{\Xi}_n(\xbf)$, initialized by $\boldsymbol{\Xi}_0=0$ is given by 
\begin{equation}
  \begin{split}
    \boldsymbol{\Xi}_n (\xbf) &= \lim_{\eps \to 0}\,\frac{\ybf_n-\xbf_n}{ \eps} \\
    %\lim_{\eps \to 0}\,{M(\ybf_{n-1},n-1) - M(\xbf_{n-1},n-1) \over \eps}=
     &= DM(\xbf_{n-1},n-1)\,  \times   \lim_{\eps \to 0} \frac{\ybf_{n-1}-\xbf_{n-1}}{ \eps} \\
     &= DM(\xbf_{n-1},\,n-1)  \boldsymbol{\Xi}_{n-1} \\
     &= \Lop_n(\xbf)\,\boldsymbol{\xi} \,.
  \end{split}
  \label{eq2-2}
\end{equation}
The random vector $\boldsymbol{\Xi}_n$ has zero mean and covariance matrix
\begin{equation}
    \Sigma^2_n(\xbf)= \mean{\boldsymbol{\Xi}_n(\xbf)\boldsymbol{\Xi}_n^T(\xbf)}= \Lop_n(\xbf)\Lop_n^T(\xbf) \,.
    \label{eq2-3}
\end{equation}
Oseledets theorem~\cite{Oseledets1961} states that the limit
\begin{equation}
    \lim_{n\to\infty}(\Lop_n^T\Lop_n)^{1/2n}= \Wop\,e^{\Lambda}\,\Wop^T 
\end{equation}
exists, where $\Wop$ is an orthogonal symplectic matrix and $\Lambda$ is diagonal with entries $\lambda_j(\xbf)$ ordered in a decreasing sequence in $j$.

The diagonal entries of $\Lambda$ are the Lyapunov exponents, and the columns of $\Wop$ the corresponding Lyapunov vectors. Since the eigenvalues of $\Lop_n^T\Lop_n$ are the same as those of the covariance matrix $\Lop_n\Lop_n^T$, the two matrices have the same characteristic polynomial. Then consider the corresponding invariants $I_n^{(k)}(\xbf), k=1,\ldots,2d$, i.e., the coefficients of the characteristic polynomial. The first invariant $I^{(1)}_n(\xbf)$, is given by the trace of the covariance matrix, namely,
\begin{equation}
  I^{(1)}_n(\xbf)\equiv \Tr\bigl(\Sigma^2_n(\xbf)\bigr)= \Tr\bigl(\Lop_n^T(\xbf)\,\Lop_n(\xbf)\bigr)=LE_n^2(\xbf) \,,
   \label{eq2-5}
\end{equation}
which is the square of the Lyapunov error $LE_n(\xbf)$. Note that it does not depend on the initial deviation vector or on the chosen orthogonal reference frame, and its asymptotic behavior is determined by the first, i.e., largest, Lyapunov exponent $\lambda_1$.

The other invariants $I^{(k)}$ are the sum of all products that combine $k$ distinct eigenvalues if they are simple. The geometric interpretation is straightforward. Letting $\ebf_j$ be the standard base vectors, we have $\Lop_n=(\ebf_{1\,\,n},\ldots,\ebf_{2d\,\,n})$where $\ebf_{j\,\,n}=\Lop_n\,\ebf_j$. %  is   the image of the base vector $j$  by the tangent map.
As a consequence, the invariant $I_n^{(k)}(\xbf)$ is the sum of the squared volumes of the $\genfrac(){0pt}{2}{2d}{k}$ parallelotopes whose sides are the vectors $\ebf_{j_1\,\,n}(\xbf),\ldots,\ebf_{j_k\,\,n}(\xbf)$.

The difference with respect to $GALI^{\,(k)}_{\,n}$ indicators (see Subsection~\ref{subsec:other}), is that the $I^{(k)}_n(\xbf)$ are  independent of the initial displacements. 

For a symplectic map, $\Lop_n(\xbf)$ is a symplectic matrix, and $\Sigma^2_n(\xbf)=\Lop_n(\xbf)\,\Lop_n^T(\xbf)$ is symplectic and positive definite. As a consequence, ordering the eigenvalues in a decreasing sequence, we have $e^{\lambda_{j;n}}e^{\lambda_{2d-j+1;n}}=1$. The asymptotic behavior of the invariant $I_n^{(k)}, k\le 2d$ is given by 
\begin{equation}
  \lim_{n\to \infty}\,\frac{1}{2n}\,\log I^{(k)}_n(\xbf) = \lambda_1(\xbf)+\ldots+\lambda_k(\xbf) \, .
  \label{eq2-7}
\end{equation}

In a region of chaotic motion, $\lambda_{j\,\,n}(\xbf)$ are positive for $j\le d$ just as their limit $\lambda_j(\xbf)$,  so that $I_n^{(k)}(\xbf)$ has exponential growth with $n$, for $n$ sufficiently large. In a region of regular motion, $I^{(k)}_n(\xbf)$ grows according to a power law $I^{(k)}_n(\xbf) \sim n^{2k}$ for $k\le d$ as all Lyapunov exponents vanish.
\subsection{Fast Lyapunov Indicator and Weighted Birkhoff averaging\label{subsec:fli}}
The Fast Lyapunov Indicator~\cite{Froeschle1997}, is one of the best known dynamic indicators, due to its straightforward implementation and its sensitiveness to the detection of chaotic structures~\cite{Lega2016fli}. Given $M(\vb{x}, n)$, its tangent map $\Lop_n(\xbf)$, and an arbitrary initial unitary deviation vector $\boldsymbol{\xi}$, $FLI$ is defined for $n\geq1$, as:
\begin{equation}
    FLI_n(\vb{x}_0, \boldsymbol{\xi}) = \ln{\norm{\Lop_n(\xbf_0) \boldsymbol{\xi}}} \, ,
    \label{eq:fli}
\end{equation}
i.e., the logarithm of the linear response $\boldsymbol{\Xi}_n(\xbf)$, calculated for an arbitrary fixed deviation vector. The quantity $FLI_n/n$ tends to the largest Lyapunov exponent as $n\to \infty$. Therefore, in a region of regular motion, this quantity tends to zero, whereas in a region of chaotic motion it takes a positive value.

It is possible to take advantage of the properties of the logarithm in Eq.~\eqref{eq:fli} to avoid overflows for large values of $n$, and to express the limit $FLI_n/n$ as an average along the trajectory $\xbf_n$~\cite{Alligood1996}:
\begin{equation}
\begin{aligned}
    \frac{FLI_n(\vb{x}_0, \boldsymbol{\xi})}{n} &= \sum_{i=0}^{n-1} \frac{\ln{\norm{\vb{y}_i-\vb{x}_i}}}{n},\\ 
    \vb{y}_i &= DM(\vb{x}_{i-1}, i-1)\frac{\vb{y}_{i-1}-\vb{x}_{i-1}}{\norm{\vb{y}_{i-1}-\vb{x}_{i-1}}} \,,\\
    \vb{y}_1 &= DM(\vb{x}_0, 0)\boldsymbol{\xi}\,.
    \label{eq:fli_mean}
\end{aligned}
\end{equation}
In the work of Das et al.~\cite{Das_2017}, it is presented how the application of the Weighted Birkhoff averaging method $\mathrm{WB}_n$~\cite{Das_2018} in the evaluation of $FLI$ can lead to superconvergence properties when applied to oscillating time series. Instead of considering equal weighting $(1/n)$, the Weighted Birkhoff averaging method uses a weighting function $w\left(\frac{i}{n}\right)$, which acts similarly to a window function in spectral analysis. A function $w(t)$ that proved to be very effective in improving the convergence of quasiperiodic time series averages~\cite{Das_2018} reads as follows:
\begin{equation}
    w(t):= 
    \begin{cases}
        \exp \left[-\frac{1}{t(1-t)}\right], & \text { for } t \in(0,1) \\ 
        0, & \text { for } t \notin(0,1)
    \end{cases} \,.
    \label{eq:birkhoff}
\end{equation}
Replacing the standard mean with $w(t)$ in Eq.~\eqref{eq:fli_mean} leads to the weighted Fast Lyapunov Indicator $FLI_n^{WB}$:
\begin{equation}
    FLI_n^{WB}(\vb{x}_0, \boldsymbol{\xi}) = \sum_{i=0}^{n-1} w\left(\frac{i}{n}\right) \ln{\norm{\vb{y}_i - \vb{x}_i}} \, .
    \label{eq:fli_birkhoff}
\end{equation}
We expect that $FLI_n^{WB}(\vb{x}_0, \boldsymbol{\xi})$ converges faster than $FLI_n(\vb{x}_0, \boldsymbol{\xi})/n$ to their common limit at least in the case of regular orbits.

To simplify the notation, in the numerical analysis we refer to $FLI_n(\vb{x}_0, \boldsymbol{\xi})/n$ and $FLI_n^{WB}(\vb{x}_0, \boldsymbol{\xi})$ as $FLI_n(\boldsymbol{\xi})/n$ and $FLI_n^{WB}(\boldsymbol{\xi})$, respectively, specifying the choice made for the initial unitary displacement $\boldsymbol{\xi}$. 
\subsection{Backward-Forward reversibility error\label{subsection:bf}}
The reversibility error is obtained by computing the linear response of the dynamics to small additive stochastic perturbations on the orbit after $n$ forward iterations $n$ followed by $n$ backward iterations 
\begin{equation}
 \begin{split}
   &\ybf_{n'}  = M(\ybf_{n'-1},n'-1) +\eps\boldsymbol{\xi}_{n'} \, ,  \quad \ybf_0=\xbf \\
   & \hskip 5 truecm 1\le n'\le n \, ; \\ 
   &\ybf_{n'} = M^{-1}(\ybf_{ n'-1}, 2n-n')+\eps   \boldsymbol{\xi}_{n' } \, , \\
   &   \hskip 5 truecm n+1\le n'\le 2n \, .
  \end{split}
\label{eq2-8}
\end{equation}
where $\boldsymbol{\xi}_{n'}$ are random vectors with zero mean and unit covariance matrix $\langle\boldsymbol{\xi}_{n'}\rangle=0$  and $\langle\boldsymbol{\xi}_{n'} \boldsymbol{\xi}_{n''}^T\rangle=\delta_{n' n''}$.
%Using the numerical finite precision as additive noise is a natural choice, but one can use an additional noise. Random displacement in backward iterations can be chosen equal to zero by setting $\boldsymbol{\xi}_{n'}=0$ for $n'>n$. 
We denote by $\xbf_{n'}$ the orbit when random deviations are absent $\eps=0$. This orbit
enjoys the symmetry property $\xbf_{n'}= \xbf_{2n-n'}$ for $n+1\le n'\le 2n$, so the reversibility condition $\xbf_{2n}=\xbf$ is satisfied.

The linear response for the $BF$ process is defined by 
\begin{equation}
   \boldsymbol{\Xi}^{BF}_{n'}(\xbf)  = \lim_{\eps\to 0}\frac{\ybf_{n'}-\xbf_{n'}}{ \eps},  \quad \qquad  1\le n'\le 2n,  %\\ \\
\label{eq2-9}
\end{equation}  
and the cumulative random deviation $\boldsymbol{\Xi}^{BF}_{n'}(\xbf)$ satisfies the recurrence 
\begin{equation}
  \begin{split}
  & \boldsymbol{\Xi}_{n'}^{BF} = DM(\xbf_{n'-1},n'-1)\,\boldsymbol{\Xi}_{n'-1}^{BF}(\xbf)+ \boldsymbol{\xi}_{n'} \, , \\ 
  & \hskip 5 truecm 1\le n'\le n;\\ 
    & \boldsymbol{\Xi}_{n'}^{BF} = DM^{-1}(\xbf_{2n-n' +1},2n-n')\, \boldsymbol{\Xi}_{n'-1}^{BF}    + \boldsymbol{\xi}_{n'} \, , \\
    & \hskip 5 truecm n+1 \le n'\le 2n \, .
  \end{split}
\label{eq2-10}
\end{equation}

From the recurrence relation of the tangent map~\eqref{eq2-1} evaluated for $n'$ and from the equality $DM^{-1}(M(\xbf,k),k))\,DM(\xbf,k)=\Iop$ for $k=2n-n'$ it follows 
\begin{equation}
  \begin{split}
  DM(\xbf_{n'-1},\,n'-1) &= \Lop_{n'}(\xbf)\,\Lop_{n'-1}^{-1}(\xbf) \, , \\
  DM^{-1}(\xbf_{2n-n' +1},2n-n') &= \Bigl(DM(\xbf_{2n-n'},2n-n')\,\Bigr)^{-1} \\
  &= \Lop_{2n-n'}(\xbf)\,\Lop^{-1}_{2n-n'+1}(\xbf) \, .
  \end{split}
\label{eq2-11}
\end{equation}
Replacing Eq.~\eqref{eq2-11} in Eq.~\eqref{eq2-8} we obtain the final result
\begin{equation}
  \begin{split}
    \boldsymbol{\Xi}_{n}^{BF}(\xbf) & = \,\Lop_n(\xbf)\,\sum_{k=1}^n\,\Lop_k^{-1}(\xbf) \,\boldsymbol{\xi}_k \,; \\  
     \boldsymbol{\Xi}_{2n}^{BF}(\xbf) &= \Lop_n^{-1}(\xbf)\,\boldsymbol{\Xi}_{n}^{BF}(\xbf)+ \sum_{k=0}^{n-1}\,\Lop^{-1}_k(\xbf)\,
    \boldsymbol{\xi}_{2n-k} \\ % = \\
    &  =\sum_{k=1}^{n-1}\,\Lop_k^{-1}(\xbf)\,(\boldsymbol{\xi}_k+\boldsymbol{\xi}_{2n-k}) \,+
   \,\boldsymbol{\xi}_{2n}+ \Lop_n^{-1}(\xbf) \boldsymbol{\xi}_n .
   \end{split}
\label{eq2-12}
\end{equation}
If random deviations are present only in the forward process, the covariance matrix of $\boldsymbol{\Xi}^{BF}_{2n}$ is given by 
\begin{equation}
  \begin{split}
    \Sigma^{2\,\,BF}_n(\xbf) &= \parmean{  \boldsymbol{\Xi}^{BF}_{2n}(\xbf)\bigl(\boldsymbol{\Xi}^{BF}_{2n}(\xbf)\bigr)^T} \\
    &=\sum_{k=1}^n \,\bigl( \,\Lop_k^T(\xbf)\Lop_k(\xbf)\,\bigr)^{-1}\, .
   \end{split}
\label{eq2-13}
\end{equation}
If random deviations are present both in the forward and backward processes, we define $\Sigma^{2\,\,BF}_n$ as $1/2$ the covariance matrix of $\boldsymbol{\Xi}^{BF}_{2n}$, and the result is the r.h.s. of Eq.~\eqref{eq2-13} where the last term of the sum $(\Lop_n^T\Lop_n)^{-1}$ is replaced by $\frac{1}{2}\, \Iop+\frac{1}{2} (\Lop_n^T\Lop_n)^{-1}$ due to the boundary condition, and asymptotically, for $n\to \infty$ the difference is negligible. 
\begin{comment}
Finally, to get another indicator, one can reverse the previous process by considering first the backward evolution defined by the map $M^{-1}(\xbf,-k)$  followed by the $F$ evolution (see the Appendix for more details). 
\end{comment}

The invariants of the matrix $\Sigma^{2\,\,BF}_n$, i.e., the coefficients of the characteristic polynomial $\det \left ( \Sigma^{2\,\,BF}_n - \lambda \Iop\right )\, , \lambda \in \mathbb{C}$ and $\lambda_i$ an eigenvalue, provide
information on the effect of small random perturbations along the orbits. If the map $M$ is symplectic, both $\Lop_n$ and $\Lop_n^T\Lop_n$ are symplectic matrices and the trace of $\Lop_n^T\Lop_n$ and its inverse are equal. As a consequence, it is not difficult to check that the trace of $\Bigl(\sum_{n'}\,\bigl(\Lop_{n'}^T\,\Lop_{n'}\bigr)^{-1}\,\Bigr)^k$ and of $\Bigl(\sum_{n'}\, \Lop_{n'}^T\,\Lop_{n'}\,\Bigr)^k$ are equal and the invariants of the covariance matrices of the $BF$ process become
\begin{equation}
 \begin{aligned}
   I^{(k)\,\,BF}_n(\xbf)&= I^{(k)} \parton{\sum_{k'=1}^n\,\Bigl(\Lop_{n'}^T(\xbf)\,\Lop_{n'}(\xbf)\,\Bigr)^{-1}} \\% = \\
   & =I^{(k)} \parton{\sum_{k'=1}^n\,\Lop_{n'}^T(\xbf)\,\Lop_{n'}(\xbf)\,} \, . \\
%    I^{(k)\,\,FB}_n(\xbf)& = I^{(k)} \parton{\sum_{k'=1}^n\,\Bigl(\Lop_{-n'}^T(\xbf)\,\Lop_{-n'}(\xbf)\,\Bigr)^{-1}} \\% = \\
%  & =  I^{(k)} \parton{\sum_{k'=1}^n\,\Lop_{-n'}^T(\xbf)\,\Lop_{-n'}(\xbf)\,} .
 \end{aligned}
\label{eq2_18}
\end{equation}

The first invariant has a very simple relation to the Lyapunov error $LE_n(\xbf)$. Explicitly, we have the following
\begin{equation}
 \begin{aligned}
   \Bigl(E^{BF}_n (\xbf))\Bigr)^2  \equiv I^{(1)\,\,BF}_n(\xbf) &= \sum_{n'=1}^n \,\Tr \Bigl(\Lop^T_{n'}(\xbf)
   \Lop_{n'}(\xbf) \Bigr) \\
   &=  \sum_{n'=1}^n \,  \Bigl(LE_n (\xbf))\Bigr)^2 \,. \\ 
%   \Bigl(E^{FB}_n (\xbf))\Bigr)^2 \equiv I^{(1)\,\,FB}_n(\xbf) &= \sum_{n'=1}^n \,\Tr \Bigl(\Lop^T_{-n'}(\xbf)
%  \Lop_{-n'}(\xbf) \Bigr) \\
%   &=  \sum_{n'=1}^n \,  \Bigl(LE_{-n} (\xbf))\Bigr)^2 \,.
 \end{aligned}
\label{eq2_19}
\end{equation}
%
\begin{comment}
For an integrable map with $d=1$, given by  $M(\xbf)= \Rop(\Omega)\xbf$  with
$\Omega=\Omega(\frac{1}{ 2}  \xbf\cdot \xbf)$
one has $\Tr(\Lop_n^T\Lop_n)= 2 +\left( \xbf\cdot\xbf\,\,\Omega'\,\,n\right)^2$ and in this case
$LE_{-n}=LE_n$ and $E^{BF}_n =LE_{-n}$. For a generic orbit of a nonintegrable system,
$LE_n$ and $LE_{-n}$ differ, although asymptotically the difference is usually negligible. 

Adding a weak dissipation, the fixed and periodic elliptic points become attractive, and the Lyapunov errors
$LE_n$ and $LE_{-n}$ and the $FB$ and $BF$ reversibility errors are substantially different,
and this is also the case for higher-order invariants.
\end{comment}

We conclude by observing that the $BF$ reversibility error analysis can be applied to investigate the effect of rounding errors in numerical computations~\cite{PANICHI201653}. Letting $M_\eps$ be the map evaluated with round-off errors and $M_\eps^{-1}$ its inverse, we have $M^{-1}_\eps(M_\eps(\xbf))= \xbf+ O(\eps)$. In the IEEE 754 international standard, the precision of a real number is $\eps\sim 10^{-16}$. Iteration with rounding is defined by Eq.~\eqref{eq2-8} where $\eps\xi_{n'}$ is missing, but $M$ is replaced by $M_\eps$. The matrix $\frac{1}{ 2}\, \boldsymbol{\Xi}_{2n}^{BF}\,\bigl(\boldsymbol{\Xi}_{2n}^{BF}\bigr)^T$, whose average defines the covariance matrix of the $BF$ reversibility error, is replaced by 
\begin{equation}
  \Xop_{2n}^{BF}(\xbf)= \frac{1}{ 2}\,\frac{ \ybf_{2n}-\xbf}{ \eps}\,\frac{(\ybf_{2n}-\xbf)^T}{\,\eps} \,.
\label{eq2_20}
\end{equation}  
This matrix has a nonzero eigenvalue, with eigenvector $\ybf_{2n}-\xbf$, and a null eigenvalue of multiplicity $2d-1$ with eigenspace orthogonal to $\ybf_{2n}-\xbf$. The noise-induced Reversibility Error Method ($REM$) squared is the nonzero eigenvalue of such a matrix, equal to its trace, and given by  
\begin{equation}
    \begin{split}
    \Bigl(REM_n^{BF}(\xbf)\Bigr)^2 &= \Tr\Bigl(   \Xop_{2n}^{BF}(\xbf)\Bigr) \\
    &= \frac{1}{ 2} \,\, \frac{\ybf_{2n}-\xbf}{ \eps} \cdot \frac{\ybf_{2n}-\xbf}{ \eps} \,.
    \end{split}
\label{eq2_21}
\end{equation}
%
%  which is the analogue of ${1\over 2}\, \boldsymbol{\Xi}\cdot\boldsymbol{\Xi}$ before taking the average, which leads to
%   $\bigl(E_{n}^{BF}(\xbf)\bigr)$.

The main difference is that $REM$, due to rounding, is the result of a single realization with a pseudorandom error and, therefore, is affected by large fluctuations when we vary $n$ or $\xbf$. These fluctuations are absent for the $BF$ reversibility error previously defined, since averaging over the random deviations is carried out. The other relevant difference is that the higher-order $REM$ invariants are zero.
\begin{comment}
We can define the indicator $REM$ for the FB process in the same way starting from the recurrence~\eqref{eq2-14} where $\boldsymbol{\xi}_{n'}$ is set to 0 and $M$ is replaced by $M_\eps$. The corresponding error is denoted by $REM_n^{FB}(\xbf)$.   %The major difference between $E^{BF}_n(\xbf)$ and
  % $REM_n^{BF}(\xbf)$ is that the first one is the result of an averaging process over the random
  % deviations whereas the second is the result of a  unique realization of a pseudo-random process.
In this case, $REM_n^{FB}(\xbf)$ is also close to $E^{FB}_n(\xbf)$ if fluctuations are washed out, for example, with a moving average in $n$. 
\end{comment}

Note that the implementation of $REM$ is trivial since it does not require the evaluation of the tangent map and the computational cost is just twice the cost of the orbit computation, provided that the inverse map is explicitly known. %This method has been used successfully in celestial mechanics. 
\subsection{$GALI^{(k)}$ indicators\label{subsec:other}}
%
  %Among the variational indicators  the fast Lyapunov indicator  was  first introduced by Froechelet et al.,
  %being  defined as
  %$FLI_{\,n}=\log \Vert\Lop_n(\xbf)\etabf\Vert$,  where $\etabf$ is a given unit vector.
  %The advantage of $E_n$ with respect to $FLI_{\,n}$ is its independence from the initial vector.
  %\\
The $k$-order indicators $GALI^{(k)}$ use the volumes of parallelotopes whose sides are normalized images of the $k$ linearly independent vectors $\boldsymbol{\eta}_j$ with $1\le j\le k$.
%
%\begin{equation}
%  \hskip -.2 truecm 
%  GALI^{(k)}_n ( \xbf)= \Vert \etabf_{1\,\,n}(\xbf)\wedge \cdots \wedge \etabf_{k\,\,n}(\xbf) \Vert
%  \quad \etabf_{j\,\,n }(\xbf)= {\Lop_n(\xbf) \etabf_j \over \Vert \Lop_n(\xbf) \etabf_j \Vert}
%  \label{eq2_22}
%\end{equation}
%
\begin{equation}
  GALI^{(k)}_n ( \xbf)= \left \Vert \frac{\Lop_n(\xbf) \boldsymbol{\eta}_1 }{ \Vert \Lop_n(\xbf) \boldsymbol{\eta}_1 \Vert}
    \wedge \ldots \wedge \frac{\Lop_n(\xbf) \boldsymbol{\eta}_k}{ \Vert \Lop_n(\xbf) \boldsymbol{\eta}_k\Vert} \right \Vert \,,
  \label{eq2_22}
\end{equation}
where $\wedge$ stands for the external product of two vectors. Their asymptotic behavior for chaotic orbits, whose first $d$ Lyapunov exponents are positive, is given by
\begin{equation}
    GALI^{(k)}_n \sim e^{-n\,\bigl((\lambda_1-\lambda_2) +\ldots+(\lambda_1-\lambda_k)\,\bigr)}  \,.
\end{equation}
where we assume a decreasing order for the exponents.

For regular, quasi-periodic  orbits,  whose Lyapunov exponents vanish, the $GALI^{(k)}$ indicators decay following a power law. We recall that the Lyapunov error invariants $I_n^{(k)}$ grow exponentially with a coefficient given by the sum of the first $k$ Lyapunov exponents for chaotic orbits, or according to a power law for regular orbits. %Growth, rather than decay with $n$, is the first significant difference with the indicators $GALI^{(k)}$. The second significant difference between the indicators $LE$ and $RE$ with respect to $GALI^{(k)}$ is their invariance with respect to the choice of the initial deviation vectors.
\begin{comment}
Regarding the indicators of reversibility error, the origin in the 
literature can be found in~\cite{} with respect to rounding. We first introduce the reversibility error as the linear response to additive noise~\cite{} to provide a plausible justification for the reversibility breakdown due to rounding. For Hamiltonian systems, the reversibility error is simply related to the Lyapunov error according to the following intuitive result.
\begin{equation}
    \bigl(\,E^{BF}_n(\xbf)\,\bigr)^2 =LE^2_1(\xbf) + LE^2_2(\xbf) +\ldots+LE^2_n(\xbf)\,.
\end{equation}
\end{comment}
%
\subsection{Introducing filters}
We conclude by remarking that the introduction of a filter such as $MEGNO$~\cite{Gozdziewski01, Cincotta2016} that drastically reduces the numerical oscillations of the indicator of chaos may greatly improve the efficiency of the indicator. In principle, the oscillations disappear using suitable normal coordinates for the considered systems, but their computation faces the limits and technical difficulties of perturbation theory. Referring to the phase flow that interpolates the orbits at integer times $t=n$, $MEGNO$, applied to $LE_t(\vb{x})$, it has the double-time average of $d\log LE_t(\vb{x})/d\log t$  %= 2t\,d \log E(t)/dt$   
\begin{equation}
\begin{split}
    MEGNO_n(LE(\vb{x})) & =  %\mean{\mean{  {d\log E^2(t)\over d\log t} }}= 
   \left\langle\left\langle  t\, \frac{d\log LE_t(\vb{x})}{ dt}\right \rangle \right \rangle\, \\
 \text{where} \qquad \mean{f(t)}& = \frac{1}{ t}\,\int_0^t \,f(t')\,dt' \, .
\end{split}
\label{eq3-23}
\end{equation}
If the indicator $LE_n(\vb{x})$ grows exponentially as $e^{\lambda t}$, then $MEGNO_n(LE(\vb{x}))$ increases as $\lambda t$. If $LE_n(\vb{x})$ follows the power law $t^\alpha$, then $MEGNO_n(LE(\vb{x}))$ converges to $2\alpha$.%, while it converges to $(2\alpha+1)/2$ for $E^{BF}(\vb{x})$.

%When $LE_n$ has exponential growth, $\mathrm{d}/\mathrm{d}n \log LE_n$ can converge to $\lambda$ more rapidly than $t^{-1}\,\log LE_n$ and faster convergence to $\lambda$ is expected for $(2t)^{-1}\,MEGNO_n(LE)$.
%
\section{\label{sec:numerical_implementations} Numerical implementations}
\subsection{Models}
To test the effectiveness of the proposed indicators of chaos, we consider a $4d$ polynomial symplectic map dependent on time, which is a generalization of the Hénon map~\cite{Yellow}. The origin is an elliptic fixed point, and the nonlinear terms combine fixed quadratic nonlinearities and variable cubic ones. The map reads:
%\begin{widetext}
\begin{equation}
    \begin{split}
    \left(\begin{array}{c}
    x_{n+1} \\
    p_{x, n+1} \\
    y_{n+1} \\
    p_{y, n+1}
    \end{array}\right) & =R(\omega_{x,n}, \omega_{y, n}) \times \\ 
    \times & \left(\begin{array}{c}
    x_{n} \\
    p_{x, n}+x_{n}^{2}-y_{n}^{2} + \mu \left(x_{n}^{3} - 3x_{n}y_{n}^{3}\right) \\
    y_{n} \\
    p_{y, n}-2 x_{n} y_{n} + \mu \left(y_{n}^{3} - 3y_{n} x_{n}^{3}\right)
    \end{array}\right) \,,
    \end{split}
    \label{eq:henon}
\end{equation}
%\end{widetext}
where $\mu$ represents the intensity of the cubic nonlinearity and $R$ is a $4\times 4$ rotation matrix defined as
\begin{equation}
    R(\omega_{x,n}, \omega_{y, n})=\left(\begin{array}{cc}
    R\left(\omega_{x, n}\right) & 0 \\
    0 & R\left(\omega_{y, n}\right)
    \end{array}\right) \,,
\end{equation}
with $R\left(\omega_{x, n}\right)$ and $R\left(\omega_{y, n}\right)$ being $2\times 2$ rotation matrices. In the following, we refer to the map~\eqref{eq:henon} as the $4d$ Hénon map and remark that it is often used as a reference model in applications such as accelerator physics (see, e.g., \cite{Yellow, invlog, Bazzani:2019csk}), since it represents the dynamics generated by a magnetic lattice that includes sextupole and octupole magnets~\cite{Yellow}.

Linear frequencies $\omega_{x, n}$ and $\omega_{y, n}$ are slowly modulated as a function of time $n$ according to
\begin{equation}
    \begin{aligned}
    &\omega_{x, n}=\omega_{x, 0}\left(1+\varepsilon \sum_{k=1}^{m} \varepsilon_{k} \cos \left(\Omega_{k} n\right)\right) \,, \\
    &\omega_{y, n}=\omega_{y, 0}\left(1+\varepsilon \sum_{k=1}^{m} \varepsilon_{k} \cos \left(\Omega_{k} n\right)\right) \,,
    \end{aligned}
\end{equation}
where $\varepsilon$ represents the modulation amplitude and the parameters $\varepsilon_k$ and $\Omega_k$ are taken from Table~1 in~\cite{invlog} to model the effect of frequency modulation in a particle accelerator due to ripples in the currents of the power supplies that feed the magnets. Modulation of the linear frequency may cause the appearance of weak chaotic regions in the stability basin near the origin. 
We recall that the parameters $\varepsilon_k$ have an order of magnitude of $10^{-4}$.

%These modulation parameters correspond to the tune modulation due to the observed ripple in the quadrupoles of the SPS~\cite{}. 
In numerical simulations, two sets of frequencies $\omega_{x0}$ and $\omega_{y0}$ have been considered, namely $(0.168,\ 0.201)$, which is close to resonances of order $5$ and $6$, and $(0.28,\ 0.31)$, which are the frequencies in the transverse phase space for charged particles orbiting in the LHC at injection energy~\cite{LHCDR}.
We have analyzed the performance of chaos indicators as a function of parameters $\varepsilon$ and $\mu$, which have been varied in the intervals $[0,64]$ and $[0,1]$, respectively. Some considerations on the computational costs of implementing the various indicators of chaos in a parallel computing architecture are reported in the Appendix~\ref{app:computing}. 

Figure~\ref{fig:survival} shows some survival plots for various configurations of the $4d$ Hénon map. A set of $300\times300$ initial conditions, sampled on a uniform Cartesian grid in the $x-y$ plane, choosing $p_x = p_y = 0$, is tracked up to $n_\mathrm{max}=10^8$ turns. Grid boundaries are selected to sample a region of interest, which depends on the linear frequencies and their modulation amplitude, that contains the stability basin of the origin, more specifically $(x, y) \in [0.0,\,0.45]\times [0.0,\,0.45]$ for case $(\omega_{x0},\,\omega_{y0}) = (0.168,\,0.201)$, or $(x, y) \in [0.0,\,0.60] \times [0.0,\,0.60]$ for case $(\omega_{x0},\,\omega_{y0}) = (0.28,\,0.31)$. An initial condition is considered stable if its distance from the origin is less than a certain control radius $r_\mathrm{c}$ when $n=n_\mathrm{max}$. Otherwise, the initial condition is considered lost and its tracking is stopped, and the stability time is given by the first value $n_\mathrm{stab}$ for which $\sqrt{x^2_{n_\mathrm{stab}}+ p_{x,n_\mathrm{stab}}^2+y_{n_\mathrm{stab}}^2+p_{y,n_\mathrm{stab}}^2}\geq r_\mathrm{c}$. The choice of $r_\mathrm{c}$ is rather arbitrary (we have considered $r_\mathrm{c}=10^2$) and the dependence of the results on $r_\mathrm{c}$ is very weak since at that amplitude the dynamics of the $4d$ Hénon map is fully dominated by polynomial terms. 
\begin{figure*}[htp]
    \centering
    \includegraphics[width=0.9\textwidth]{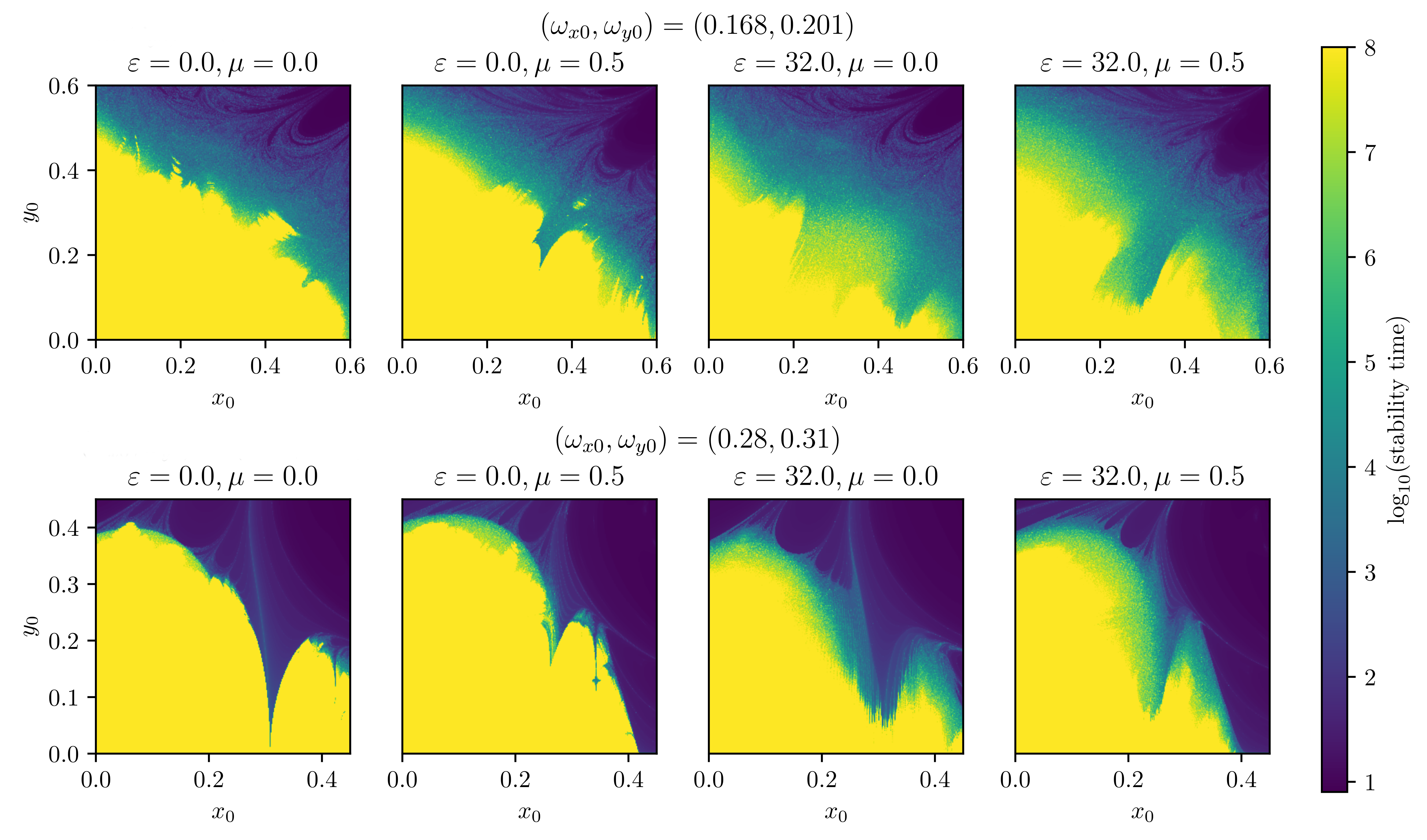}
    \caption{Survival plot for various $4d$ modulated Hénon maps with quadratic and cubic nonlinearities. Initial conditions, sampled on an uniform Cartesian $300\times300$ grid in the $x-y$ plane, are tracked up to $n_\text{max}=10^8$ and are considered lost when their distance to the origin exceeds a predefined maximum radius $r_\mathrm{c}=10^2$. The two sets of linear frequencies feature different shapes of the stable region as can be seen by comparing the plots in the two rows. The parameters $\varepsilon$ and $\mu$ induce additional changes, in particular the increase of the size of the transition region between stability up to $n_\mathrm{max}$ and shorter stability time. The color scale is related to the logarithm of the stability time as reported on the right.}
    \label{fig:survival}
\end{figure*}
The two rows of Fig.~\ref{fig:survival} show the survival plots for the two sets of frequencies considered in the studies. The shape of the stable region (yellow area) strongly depends on the frequencies, as different sets of resonances affect the dynamics. Furthermore, the impact of $\varepsilon$ and $\mu$ is also clearly seen. The first enlarges the transition region between stable initial conditions and unstable ones, i.e., the region for which $n_\mathrm{stab} < n_\mathrm{max}$, where a weak diffusion occurs, while the latter changes the shape of the stable region.
\section{\label{sec:results} Results of numerical investigations}
In the following, we report the results of the numerical study of the dynamic indicators presented in Section~\ref{sec:review}, namely $\log_{10}(LE)$, $FLI$, $FLI^{WB}$, $MEGNO(LE)$, $GALI^{(4)}$, $REM$, and $FMA$. Note that we consider the logarithm of $LE$, as it is a quantity comparable to $FLI$ and $MEGNO(LE)$. We first focus on the dependence of $FLI$ on the choice of the initial displacement vector $\vb{\xi}$, and compare it with $\log_{10}(LE)$. Next, we discuss a comparison between the convergence rate of $FLI$ and that of $FLI^{WB}$. Finally, we compare the classification performance of all dynamic indicators by determining their accuracy, together with its time dependence, in reconstructing a Ground Truth (GT) evaluated at a high iteration time.
\subsection{Dependence on the initial displacement}

\begin{figure*}[htp]
    \centering
    \includegraphics[width=\textwidth]{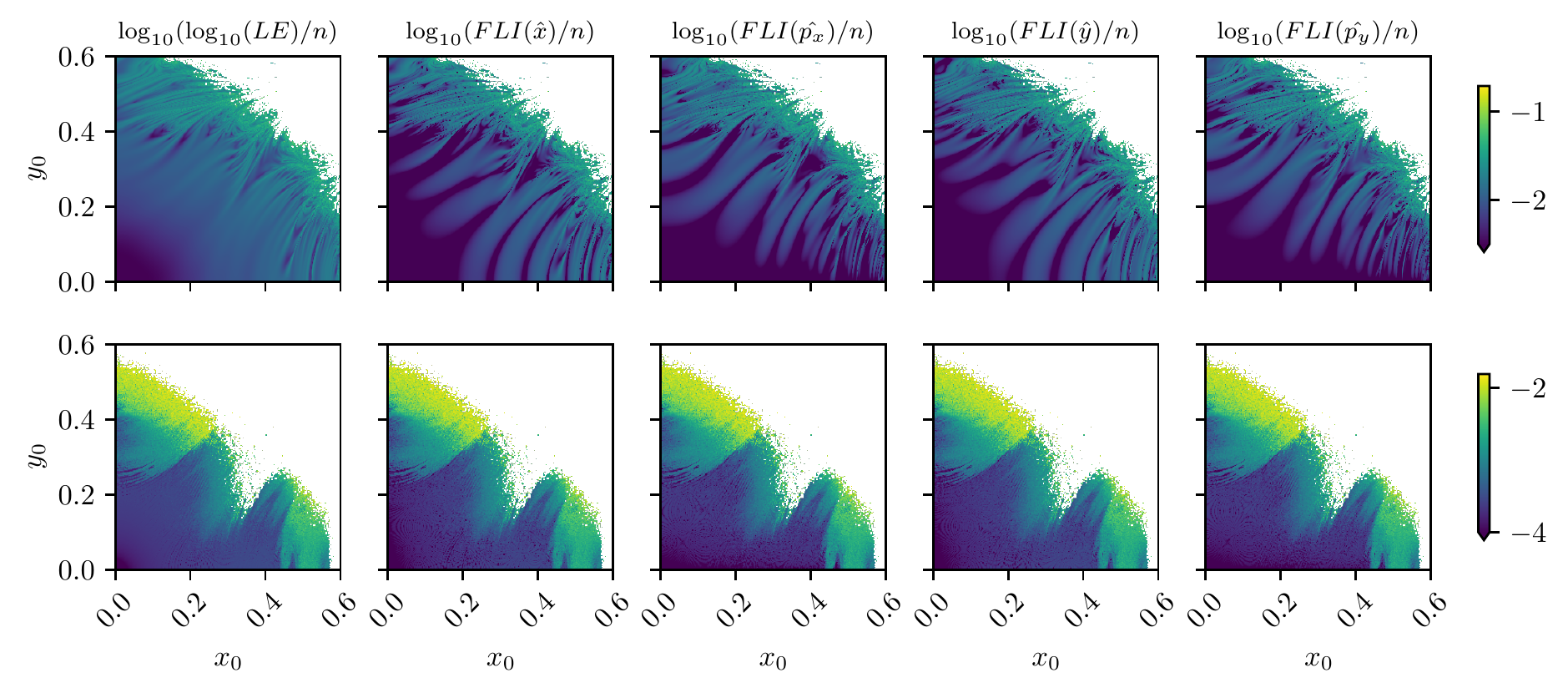}
    \caption{Color maps of $\log_{10}(\log_{10}(LE)/n)$ and $\log_{10}(FLI/n)$ indicators for a low iteration number ($n=10^2$, top row) and a high iteration number ($n=10^4$, bottom row). In both rows, the $FLI$ for the four possible displacements is shown together with $LE$, to highlight the different structures shown by the indicators. It is possible to see how, for low iteration numbers, different choices of initial displacement for $FLI$ highlight structures that do not appear in $LE$. The differences reduce for higher number of turns, but are still present. These structures have a vage resemblance with the net of resonances that is present in phase space. Note that an arrow at the bottom of the color bar means that pixels of the bottom color correspond to a value equal to or lower than the bottom value. White pixels correspond to initial conditions whose distance from the origin has exceeded a predefined radius ($r_c=10^2$) during the tracking, before reaching the target iteration number $n$. (Simulation parameters used: $(\omega_{x0},\omega_{y0})= (0.168,\ 0.201),\ \varepsilon=64.0,\ \mu=0.5$).}
    \label{fig:le_fli_compare_short_long}
\end{figure*}

The main feature of $LE$, compared to $FLI$, is its independence from the initial choice of direction of the unitary displacement vector $\vb{\xi}$. To highlight this, in Fig.~\ref{fig:le_fli_compare_short_long}, we directly compare the calculated values of $\log_{10}(\log_{10}(LE)/n)$ with those of $\log_{10}(FLI/n)$, calculated with an initial displacement along one of the four orthonormal base vectors $\hat{x},\,\hat{p}_x,\,\hat{y},\,\text{and }\hat{p}_y$. These calculations are carried out for a set of $300\times300$ initial conditions, sampled on a uniform Cartesian grid in the $x-y$ plane. It is possible to see how, at low turn number ($n=10^2$, top row), the different choice of displacement highlights the structures in $FLI$ that are missing in $LE$. This can be explained by considering that the displacement vector is not fully aligned along the largest Lyapunov exponent yet. In contrast, these structures are missing for $LE$, which has smoother behavior.

The observed differences are much reduced for a higher number of turns ($n=10^4$, bottom row), as the initial displacement tends to become almost aligned along the direction corresponding to the largest Lyapunov exponent. However, despite the smaller differences between $\log_{10}(\log_{10}(LE)/n)$ and $\log_{10}(FLI/n)$, the behavior of the various indicators is still not the same. It is worth noting how displacements along $\hat{x}$ and $\hat{y}$ produce similar structures that are, however, different with respect to the case in which displacement is carried out along $\hat{p}_x$ or $\hat{p}_y$. Globally, these observations underline the value of the invariance properties of $LE$, which seems to be more promising than $FLI$ for the analyses that will be discussed in the following sections.

As this dependence on the initial displacement decreases with higher iteration numbers, we will focus only on $FLI(\hat{x})$ for the remainder of the paper, as the rest of the results are not significantly affected by this choice.
\subsection{Application of Weighted Birkhoff averaging to $FLI$}
\begin{figure}[htp]
    \centering
    \includegraphics[width=0.5\columnwidth]{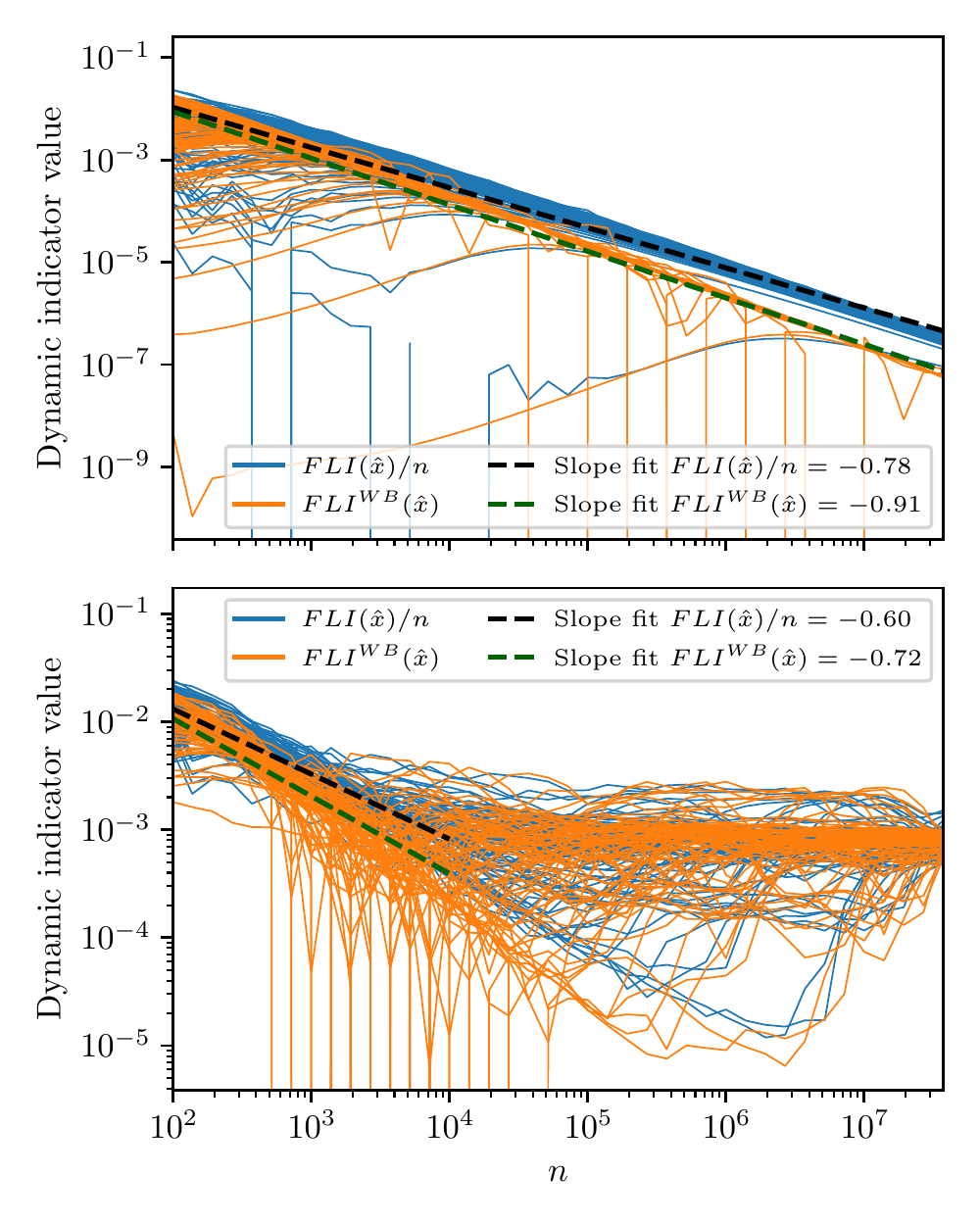}
    \caption{Time evolution of $FLI$ computed using either a standard mean or the Birkhoff averaging. Top plot: indicators computed for a set of 100 regular initial conditions, the fit highlights a faster convergence rate for the Birkhoff averaging. Bottom plot: indicators computed for a set of 100 chaotic initial conditions. A similar improvement in convergence rate is observed for low $n$ values, before reaching a saturation value of the indicator of the order of $10^{-3}$. (Simulation parameters: $(\omega_{x0},\omega_{y0})= (0.28,\ 0.31),\ \varepsilon=32.0,\ \mu=0.5$).}
    \label{fig:fli_compare_mean_birk}
\end{figure}

As an additional analysis of the time dependence of chaos indicators, we compare the values obtained for $FLI$ at different times, using the standard approach that considers the mean in Eq.~\eqref{eq:fli_mean}, that is, $FLI/n$, or the variant based on the use of Birkhoff weights as in Eq.~\eqref{eq:fli_birkhoff}, that is, $FLI^{WB}$. The analysis starts considering two ensembles of regular and chaotic particles that have been classified by means of the value of the $FLI$ indicator computed for $n=10^8$ turns (effectively this sets a ground-truth level, as discussed in the next section). The sets are also used to calculate the time evolution of $FLI/n$ and $FLI^{WB}$ with the objective of evaluating possible improvements in the latter compared to the first. In Fig.~\ref{fig:fli_compare_mean_birk} (top), the comparison is made for a subset of the set of regular initial conditions, whereas the behavior of chaotic ones is shown in the bottom plot of the same figure. It is possible to observe how, for regular initial conditions, Birkhoff averaging consistently speeds up the convergence of $FLI^{WB}$ to zero.

The case of chaotic initial conditions has different characteristics. In fact, a saturation region is observed for the indicator value on the order of $10^{-3}$ for both indicators. When this value is reached, both indicators oscillate around it. However, the slope with which this nonzero value is reached is different for the two indicators and is higher in absolute value for $FLI^{WB}$ than for $FLI/n$, similar to what is observed for the case of regular orbits. It is also worth stressing the presence of initial conditions that, up to some $n=10^6$ turns, feature a steady decrease in the value of the dynamic indicator, as if they were characterized by regular motion. However, after that, the value of the indicator suddenly increases, reaching the value that identifies chaotic orbits. This behavior clearly defies any approach aimed at classifying initial conditions as regular or chaotic in finite time. 

\begin{figure}[htp]
    \centering
    \includegraphics[width=0.5\columnwidth]{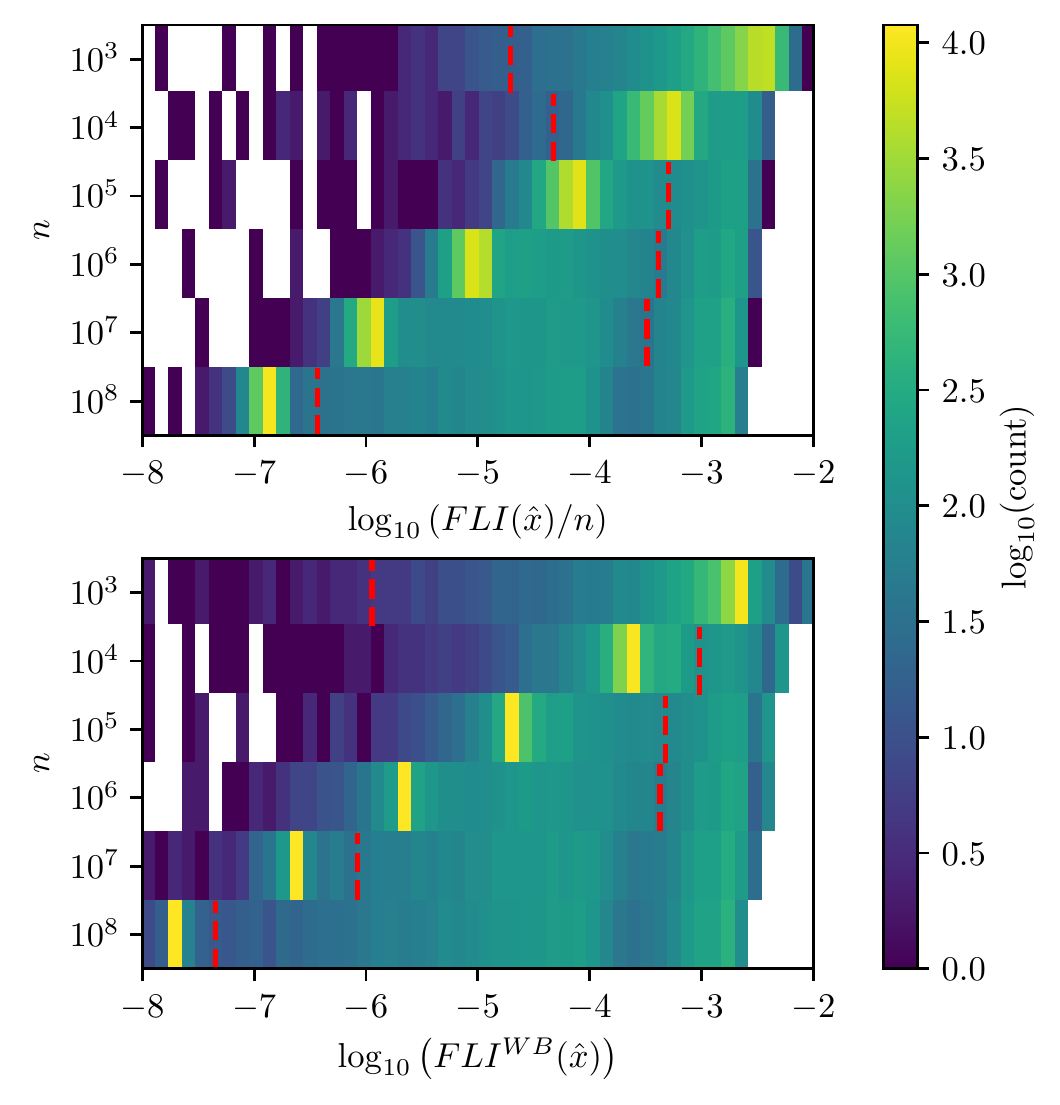}
    \caption{Time evolution of the distribution of the values of $FLI$ (top) and $FLI^{WB}$ (bottom) indicators for the whole set of 15684 initial conditions that survived up to $n_\text{max}=10^8$. The Birkhoff averaging leads to faster convergence towards zero of the regular initial conditions, which are represented by the yellow band. Furthermore, the width of such a band is narrower for $FLI^{WB}$ with respect to $FLI$. The red dashed lines represent threshold values, defined by our algorithm, representing the attempt to perform the binary classification in regular and chaotic initial conditions. (Simulation parameters: $(\omega_{x0},\omega_{y0})= (0.28,\ 0.31),\ \varepsilon=32.0,\  \mu=0.1$).}
    \label{fig:fli_colormap_mean_birk}
\end{figure}

The improvement caused by the Birkhoff averages is also clearly visible in Fig.~\ref{fig:fli_colormap_mean_birk}, where the time evolution of the distribution of the values of $FLI/n$ (top) and $FLI^{WB}$ (bottom) is shown. The part of the distribution corresponding to the regular initial conditions reaches its peak (yellow band) and moves towards zero with increasing $n$. However, the displacement towards zero is faster for $FLI^{WB}$. Furthermore, the peak of the distribution is sharper for $FLI^{WB}$ than for $FLI$. In both graphs, a faint trace of a peak is visible corresponding to the indicator value of about $10^{-3}$. This feature is remarkably similar for the two indicators, as already seen in Fig.~\ref{fig:fli_colormap_mean_birk}.

This behavior shows that the regular orbits benefit from the use of the Birkhoff averages, whereas the chaotic ones are mostly unaffected by the special averaging mechanism. These features can be exploited for the classification problem that will be addressed in the next section.

\subsection{Classification performance}\label{subsec:classification}

For this analysis, we study the predictive performance of chaos indicators in terms of a binary classification of a large set of initial conditions by varying the number of iterations $n$. It should be stressed that this classification is performed only on the orbit of an initial condition that has been detected to be stable for $n_\mathrm{max}$.

An overview of the time dependence of the dynamic indicators and the distribution of their values observed in our numerical investigation is given in the Appendix~\ref{app:timedep}. The main feature of interest, which constitutes the basis of this analysis, is the general tendency of dynamic indicators to create a bimodal distribution, as has also been reported for finite-time Lyapunov exponents in~\cite{PhysRevE.60.2761,VALLEJO200326}. We focus on studying the evolution of this specific characteristic, i.e., the presence of two peaks in the distribution of indicator values, as a function of time, which is the key feature used for the classification analysis.

As the development of the bimodal distribution requires various orders of magnitude of the number of turns, we perform our analysis on the logarithm of the seven dynamic indicators, namely $\log_{10}(\log_{10}(LE)/n)$, $\log_{10}(MEGNO(LE)/n)$, $\log_{10}(FLI/n)$, $\log_{10}(FLI^{WB})$, $\log_{10}(GALI^{(4)})$, $\log_{10}(REM)$, and $\log_{10}(FMA)$. The factor $n^{-1}$ is included in the first two indicators to observe a comparable evolution of values over time with the two $FLI$ indicators, since, ultimately, its presence does not alter the outcome of these studies.

\begin{figure*}[htp]
    \centering
    \includegraphics[width=0.9\textwidth]{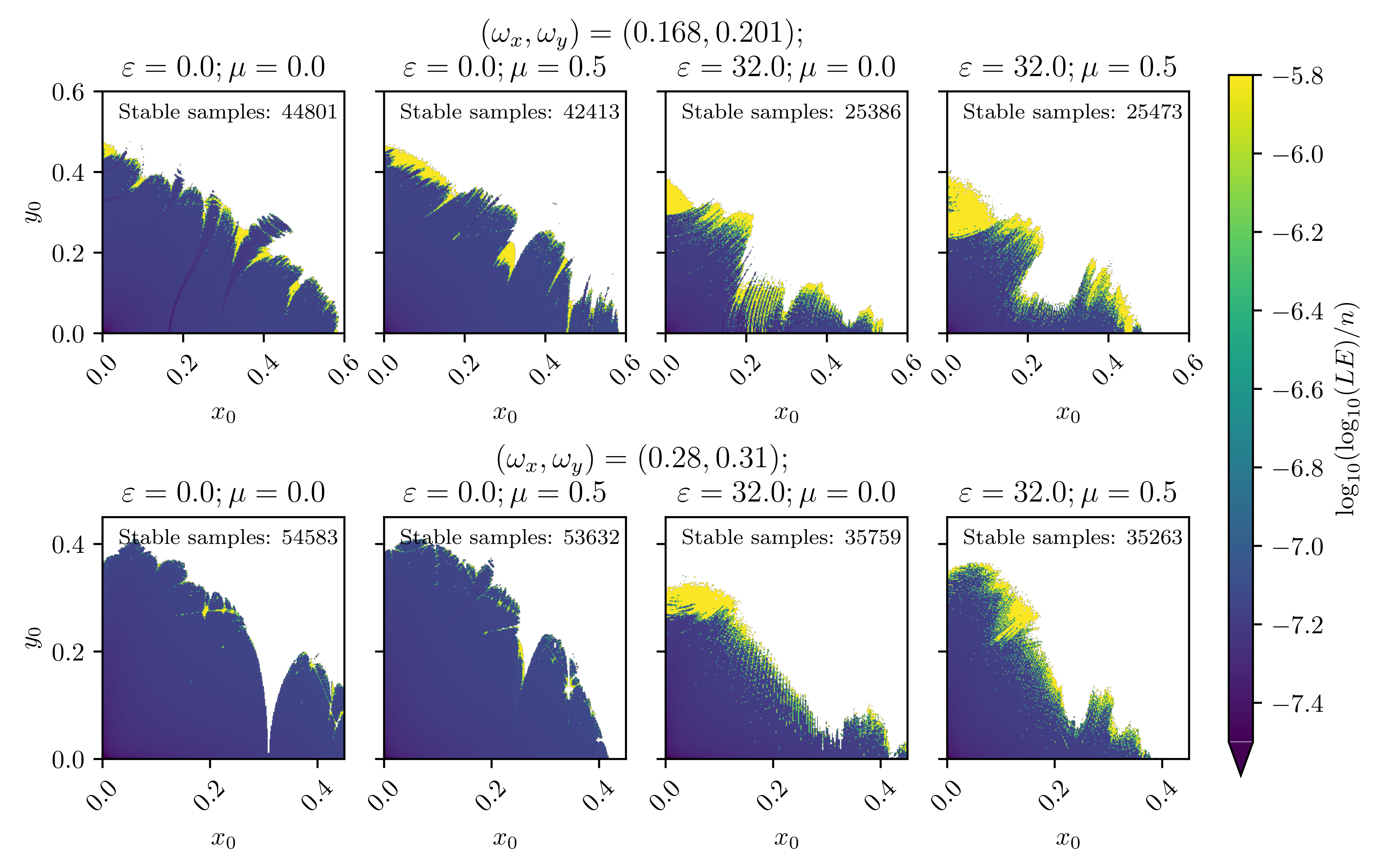}
    \caption{Distributions of $\log_{10}(\log_{10}(LE)/n)$ for various $4d$ modulated Hénon maps (the same cases shown in Fig.~\ref{fig:survival}) with quadratic and cubic nonlinearities. $300\times300$ initial conditions, sampled on an uniform Cartesian grid in the $x-y$ plane, are tracked up to $n_\text{max}=10^8$. It is possible to observe how the case for $\varepsilon=0.0,\ \mu=0.0$, corresponding to the absence of modulation and cubic nonlinearities, lead to regular motion almost everywhere, except for a small set of initial conditions. For the other cases, extended regions of chaotic motion are visible. Note that the maximum value registered in the color maps corresponds to numerical saturation.}
    \label{fig:ground_truth_bis}
\end{figure*}

To carry out this task, we first construct a ground truth (GT) for different sets of parameters for the $4d$ Hénon map, iterated for $n_\text{max} = 10^8$. The initial conditions are then classified into a binary chaotic/regular classification scheme using the $LE$ indicator. An example is given in Fig.~\ref{fig:ground_truth_bis} where eight cases, the same as those depicted in Fig.~\ref{fig:survival}, are displayed. Dark colors identify regular regions of the phase space, whereas lighter colors denote chaotic regions. It is clearly seen that the frequency modulation and the presence of the cubic nonlinearity increase the extent of the chaotic areas of the phase space, also generating regions in which regular and chaotic orbits are deeply intertwined.  

The GT classification is built from the distribution of the values of $\log_{10}(\log_{10}(LE)/n)$ for $n_\text{max}$. The resulting distribution has a main group of regular initial conditions with low value $LE$, and a second group of chaotic initial conditions with higher value $LE$. Due to the large separation of these two clusters, a threshold value has been calculated to distinguish them using a kernel density estimation method (KDE)~\cite{doi:10.1080/24709360.2017.1396742, refId0} with a Gaussian kernel and different bandwidth values. This allows investigating the Mode Tree~\cite{10.2307/1390955} of the distribution, detecting its two main modes, and setting the position of the minimum of the distribution between them. It is worth stressing that more refined approaches might be devised to detect the peaks or, equivalently, cluster the indicator values, but they have not been considered in this analysis. In fact, our focus is on the performance of the indicator in generating a suitable distribution for the classification problem, even for low values of $n$, not on designing a sophisticated algorithm to analyze the distribution of the indicator, including its peculiarities.

\begin{figure*}[htp]
    \centering
    \includegraphics[width=1.0\textwidth]{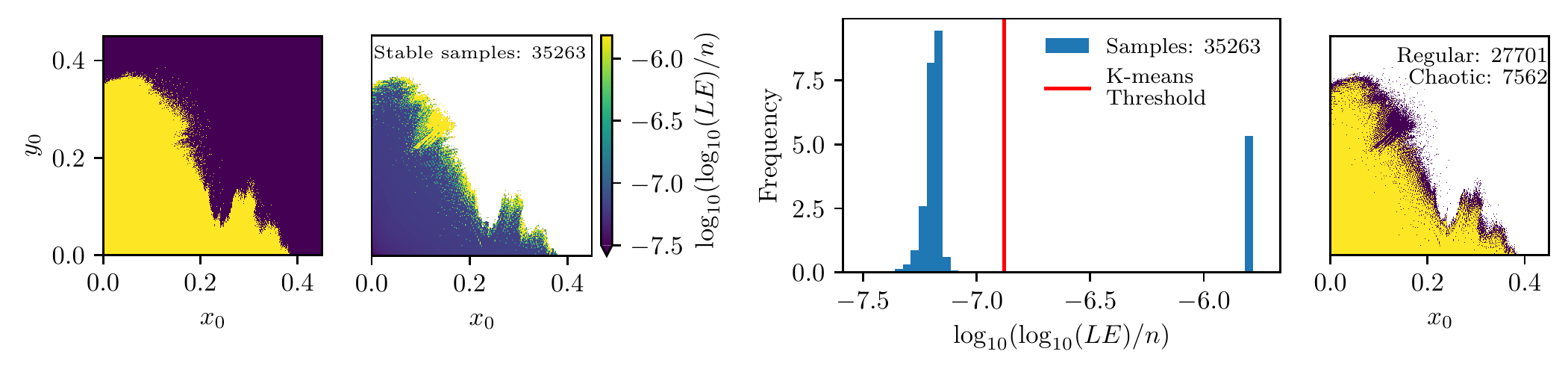}
    \caption{Ground Truth construction for a modulated Hénon map. From left to right: a survival plot of the initial conditions stable up to $n_\text{max}=10^8$ (yellow is stable, purple is unstable); distribution of the $LE$ indicator for all stable initial conditions, evaluated at $n_\text{max}$; histogram of $\log_{10}(\log_{10}(LE)/n)$, classified with a threshold evaluated with a KDE-based procedure; binary classification of regular (yellow) and chaotic (purple) initial conditions. (Simulation parameters:  $(\omega_{x0},\omega_{y0})= (0.28,\ 0.31),\ \varepsilon=32.0,\ \mu=0.5$).}
    \label{fig:ground_truth_2}
\end{figure*}

An example of the GT construction process can be seen in Fig.~\ref{fig:ground_truth_2}. 
Stable initial conditions up to $n_\mathrm{max}$ are identified by direct tracking (first graph from the left), and the value of the indicator $LE$ is calculated for the set of stable initial conditions (second graph from the left). At this stage, it is possible to compute the distribution of $LE$ and determine the threshold that separates the peaks of the bimodal distribution (third plot from the left) and provides the criterion to classify any given initial condition as regular or chaotic. Applying the computed threshold, it is possible to generate a binary map with the resulting classification (fourth plot from the left). In this case, the determination of the threshold for the case shown is rather straightforward, as the large separation between the two peaks makes the actual value of the threshold not particularly relevant. However, when $n \ll n_\mathrm{max}$ the separation between the peaks decreases and the threshold value becomes essential for an efficient classification of the initial conditions. 

\begin{figure}[htp]
    \centering
    \includegraphics[width=0.5\columnwidth]{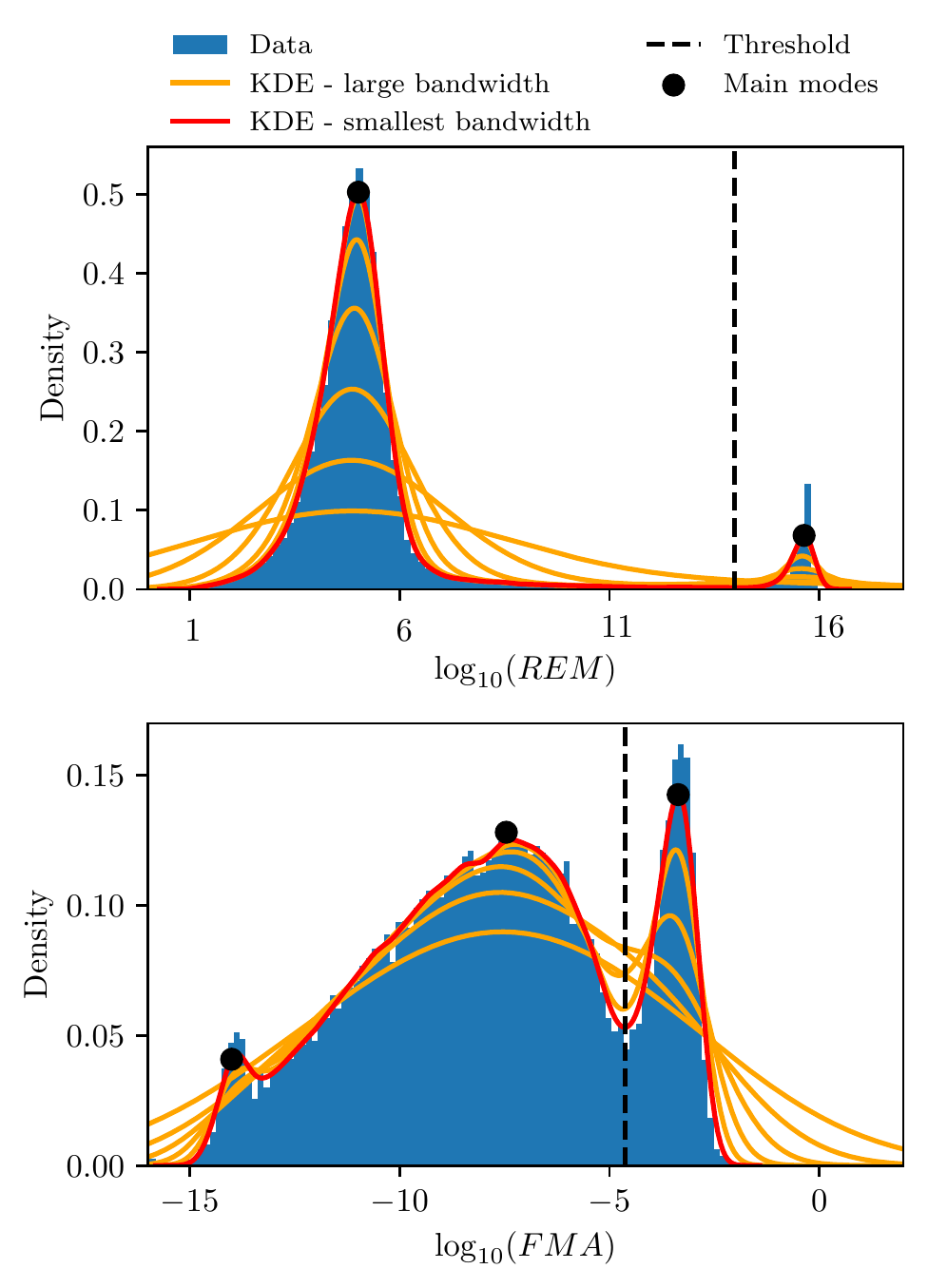}
    \caption{Example of the KDE-based procedure for computing a threshold for the binary classification (regular/chaotic) of initial conditions. Top: application to $\log_{10}(REM)$ (evaluated at $n=10^5$). KDEs with various bandwidth are used until the two main peaks of the bi-modal distribution are detected, the threshold is then placed at the position of the minimum of the distribution between them. This procedure is applied to all dynamic indicators except for $FMA$. Bottom: application of the procedure to $\log_{10}(FMA)$ (evaluated at $n=10^5$), which clearly exhibits a three-mode distribution. The procedure is applied so that it detects the three main modes of the distribution, and then sets the threshold at the minimum of the distribution between the two modes at higher values. (Simulation parameters:  $(\omega_{x0},\omega_{y0})= (0.28,\ 0.31),\ \varepsilon=32.0,\ \mu=0.5$).}
    \label{fig:thresholds}
\end{figure}

Examples of the procedure for determining the threshold based on the indicator distribution are shown in Fig.~\ref{fig:thresholds}. In the top plot, the case of $REM$ is depicted (but it is representative of all other indicators except $FMA$). The use of KDE with different bandwidth clearly shows how the two peaks of the distribution can be detected. This allows the position of the threshold to be set at the location of the minimum value of the distribution in between the two peaks. The case of $FMA$ is different since the distribution has three peaks and the standard algorithm to determine the threshold must be adapted. Therefore, KDE is used to determine the position of the three peaks, and the threshold is set at the position of the minimum of the distribution in between the two peaks with the largest amplitude.  

This choice is somewhat arbitrary, but the features of the distribution clearly indicate that the performance of the indicator is limited, with little possibility of improving it. Indeed, the non-negligible fraction of initial conditions that generate the part of the distribution in between the extreme peaks cannot be clearly classified by the proposed approach, as some of them will turn chaotic, whereas other regular if the indicator would be computed over a longer time span.  

Once the GT has been computed, we define as predictive performance of a dynamic indicator the accuracy in reconstructing the binary classification in the GT, that is, the ratio between the correctly labeled initial conditions and the total number of stable initial conditions. Such a reconstruction is attempted using the same strategy implemented for the determination of the GT, namely, we consider the distribution of the dynamic indicator under consideration and define a binary classification using a threshold computed via the KDE-based approach. The resulting thresholds evaluated over time for $REM$ and $FMA$ are visualized in detail in Fig.~\ref{fig:neo_evolution}, while the results for the other dynamic indicators are presented in Appendix~\ref{app:timedep}.

\begin{figure*}[htp]
    \centering
    \includegraphics[width=1.0\textwidth]{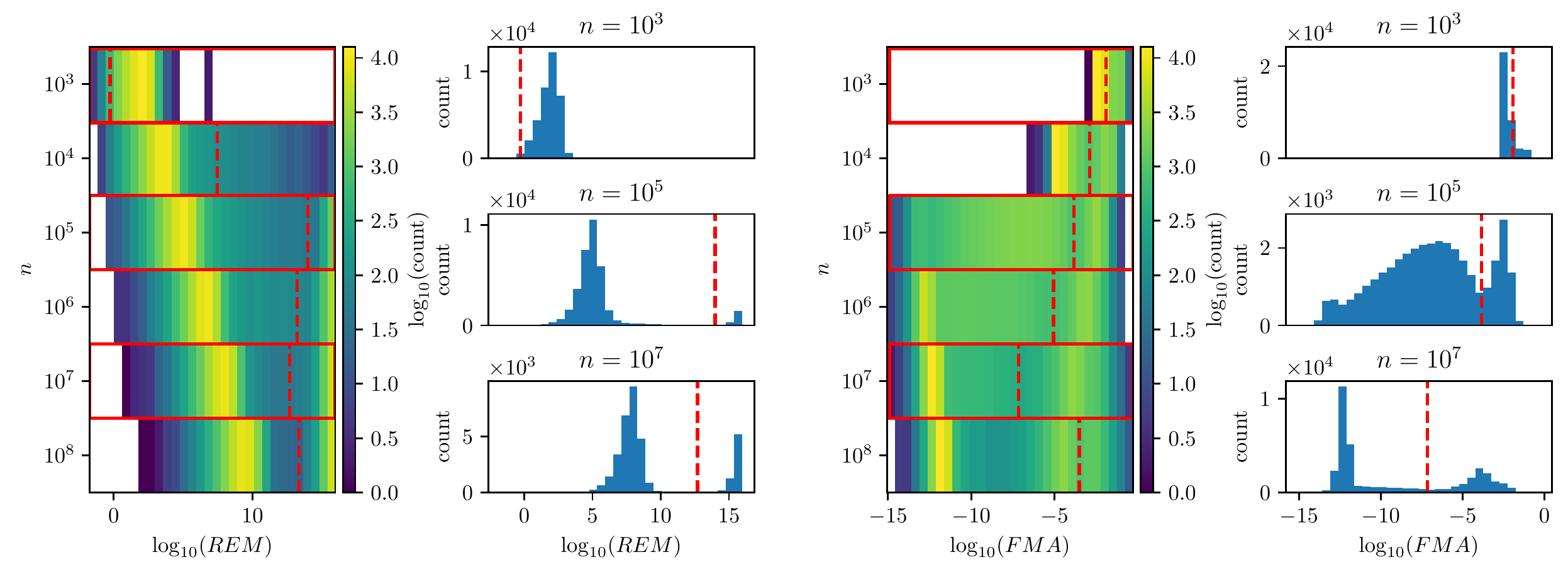}
    \caption{Distribution of values of $\log_{10}(REM)$ (left) and $\log_{10}(FMA)$ (right) as a function of time for a modulated $4d$ Hénon map. The red dashed lines represent threshold values, defined by our algorithm shown in Fig.~\ref{fig:thresholds}, representing our criterion to distinguish regular and chaotic orbits.
    For low values of the iterations $n$, the distribution of both indicators is in general represented by a uni-modal function. For higher values of $n$, we can see the formation of two separate clusters in the case of $REM$, making the distribution bi-modal. For $FMA$, we have in general a different behavior, as it tends to form a tri-modal distribution. (simulation parameters: $(\omega_{x0},\omega_{y0})= (0.28,\ 0.31),\ \varepsilon=32.0,\ \mu=0.5$).}
    \label{fig:neo_evolution}
\end{figure*}

The accuracy performance of the dynamic indicator is then evaluated for various $n < n_{\text{max}}$. We expect a good-performing dynamic indicator to achieve high accuracy values when it generates two separate groups, even when $n \ll n_{\text{max}}$. Such behavior, in fact, enables effective mode detection and consequent effective GT reconstruction. In contrast, a poor-performing dynamic indicator will need a longer tracking time before showing the presence of two separate clusters, causing the threshold determination to be unable to separate the chaotic from the regular initial conditions.

\begin{figure*}[htp]
    \centering
    \includegraphics[width=1.0\textwidth]{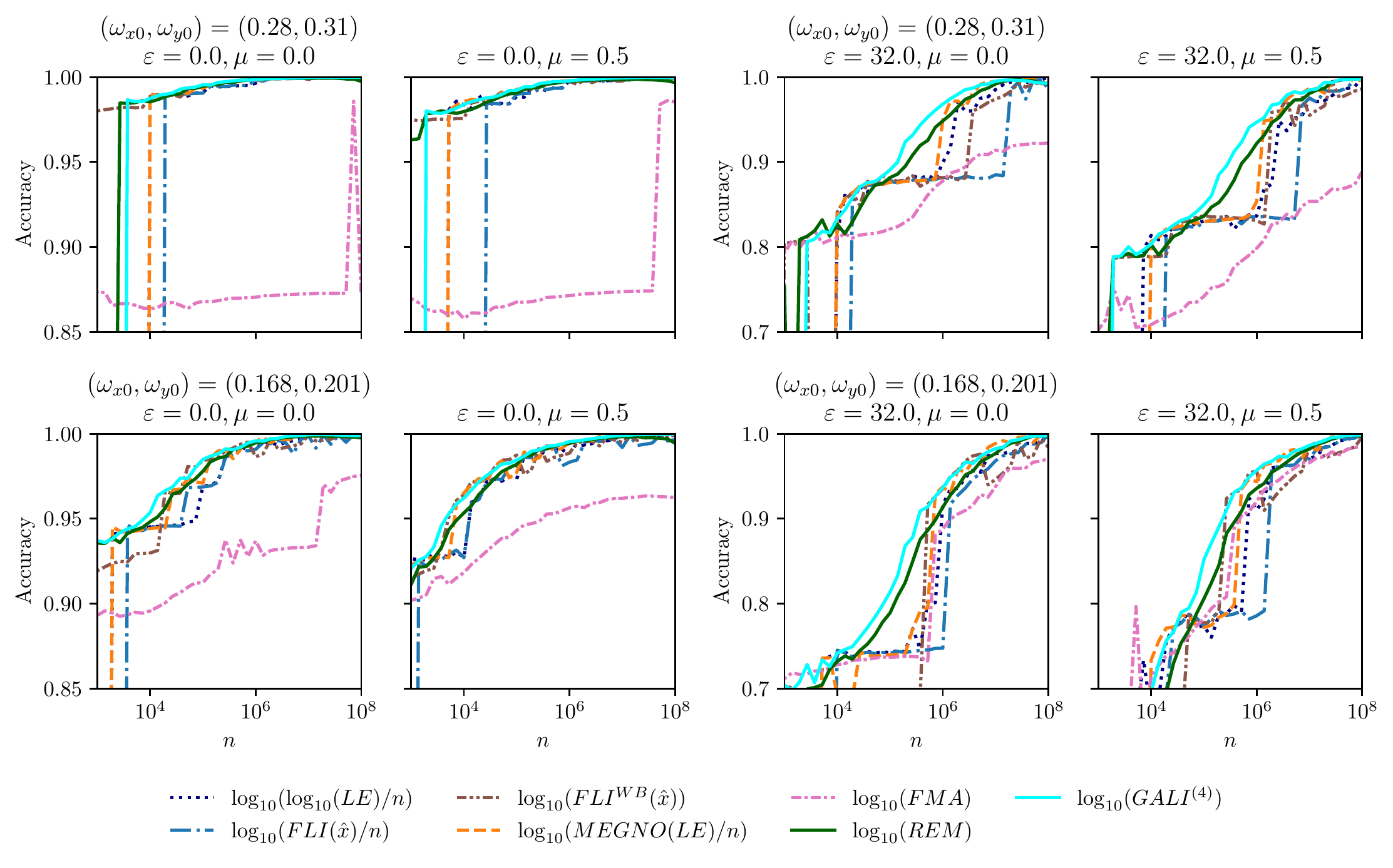}
    \caption{Time dependence of the accuracy achieved in reconstructing the Ground Truth (computed for $n_\text{max}=10^8$) by the various dynamic indicators for eight cases of the $4d$ modulated Hénon maps (the same cases shown in Fig.~\ref{fig:survival}), differing by cubic nonlinearities and frequency modulation.}
    \label{fig:performance}
\end{figure*}

A global comparison of the classification performance of the seven dynamic indicators is carried out, and the accuracy achieved by the dynamic indicators as a function of $n$ is shown in Fig.~\ref{fig:performance}, for different sets of parameter values for the $4d$ Hénon maps. 

When considering the Hénon maps with $\varepsilon=0.0$, i.e., without frequency modulation, a rather small fraction of chaotic orbits with a very mild dependence on $n$ of the accuracy of the various dynamic indicators is observed. Furthermore, $FMA$ differs from all other indicators, clearly showing poorer performance in terms of accuracy. All other indicators have very similar performance, the only difference being in the time at which a steplike increase in accuracy is observed, which occurs for $n = 10^3 - 10^4$, corresponding to 4-5 orders of magnitude lower than $n_\text{max}$. This sudden increase in accuracy is related to the time required by dynamic indicators to generate a bimodal distribution that can be efficiently analyzed using our KDE-based procedure. In this sense, it should be noted that $GALI^{(4)}$ is the most accurate indicator, as it reaches high accuracy values even at very low values of $n$ and the gradual increase does not occur in the range of $n$ shown in the graphs. In general, the behavior observed for all indicators (except $FMA$) shows that a rather accurate prediction of GT can be achieved using the information provided by the indicators over a rather limited number of turns.

In the case with $\varepsilon=32.0$, i.e., with frequency modulation and a larger fraction of chaotic orbits, the situation changes dramatically. Accuracy depends rather strongly on $n$, suggesting that chaos detection requires a larger number of turns to be accurate. In terms of the ranking of the indicators, $FMA$ remains the worst (this is certainly true for case $(\omega_{x0}, \omega_{y0})=(0.28, 0.31)$, while for case $(\omega_{x0}, \omega_{y0})=(0.168, 0.201)$ a better performance is observed). $REM$ and $GALI^{(4)}$, are the best values in a wide range of values of $n$. Furthermore, they do not show any sudden jump in accuracy because of their well-behaved distribution. Finally, we remark that beyond $n = 10^6 - 10^7$, the precision of all indicators is very similar.

To provide a quantitative assessment of the performance of the dynamic indicators, we define a performance estimate as
\begin{equation}
    \frac{1}{2}\int_4^6 \text{Accuracy}(10^x) \,\mathrm{d}x \,,
    \label{eq:final_evaluation}
\end{equation}
i.e., the integral of the accuracy achieved and displayed in Fig.~\ref{fig:performance} normalized to the integral of the ideal case with unit accuracy throughout the turn interval. The reasons for such a definition are twofold: First, it avoids the possible bias introduced by indicators that are more efficient in detecting the chaotic behavior at low number of turns but that are not so efficient afterwards; second, it probes the predictive power of the indicator by setting an upper bound that is lower than the turn number used for determining the GT. Equation~\eqref{eq:final_evaluation} has been numerically evaluated using the trapezoidal rule and considering 50 values of $n$ equally spaced on a logarithmic scale over the interval $10^4-10^6$.  The performance estimate values for the dynamic indicators for the various H\'enon maps are reported in Table~\ref{table:values}. 

\begin{table*}[htb]
    \centering
    \caption{Performance estimate of the dynamic indicators for the various Hénon map configurations, evaluated using Eq.~\eqref{eq:final_evaluation} over the interval $n=10^4-10^6$. Values are ranked in decreasing order. It is clearly seen that $GALI^{(4)}$ is the highest scorer and $REM$ is the second-best scorer in most of the cases considered. The uncertainty in the performance estimate is evaluated by applying a variation of the calculated thresholds of $\pm 5\%$.}
    \begin{tabular}{lc|lc}
        \hline
        \multicolumn{4}{c}{$(\omega_x, \omega_y) = (0.28, 0.31)$} \\
         \multicolumn{2}{c|}{$\varepsilon = 0.0; \mu = 0.0$} & \multicolumn{2}{c}{$\varepsilon = 0.0; \mu = 0.5$} \\
        \hline
        $\log_{{10}}(GALI^{{(4)}})$ & $0.99700 \pm 0.00014$ & $\log_{{10}}(GALI^{{(4)}})$ & $0.9956 \pm 0.0002$ \\ 
$\log_{{10}}(FLI^{{WB}}(\hat{{x}}))$ & $0.9966 \pm 0.0005$ & $\log_{{10}}(REM)$ & $0.99423 \pm 0.00004$ \\ 
$\log_{{10}}(MEGNO(LE)/n)$ & $0.9965 \pm 0.0008$ & $\log_{{10}}(\log_{{10}}(LE)/n)$ & $0.99 \pm 0.03$ \\ 
$\log_{{10}}(REM)$ & $0.99629 \pm 0.00003$ & $\log_{{10}}(FLI^{{WB}}(\hat{{x}}))$ & $0.99 \pm 0.10$ \\ 
$\log_{{10}}(\log_{{10}}(LE)/n)$ & $0.99 \pm 0.05$ & $\log_{{10}}(FLI(\hat{{x}})/n)$ & $0.9080 \pm 0.0016$ \\ 
$\log_{{10}}(FLI(\hat{{x}})/n)$ & $0.94 \pm 0.01$ & $\log_{{10}}(MEGNO(LE)/n)$ & $0.90 \pm 0.14$ \\ 
$\log_{{10}}(FMA)$ & $0.8738 \pm 0.0005$ & $\log_{{10}}(FMA)$ & $0.8797 \pm 0.0004$ \\ 
        \hline
        \multicolumn{4}{c}{$(\omega_x, \omega_y) = (0.28, 0.31)$} \\
         \multicolumn{2}{c|}{$\varepsilon = 32.0; \mu = 0.0$} & \multicolumn{2}{c}{$\varepsilon = 32.0; \mu = 0.5$} \\
        \hline
        $\log_{{10}}(GALI^{{(4)}})$ & $0.9453 \pm 0.0014$ & $\log_{{10}}(GALI^{{(4)}})$ & $0.924 \pm 0.002$ \\ 
$\log_{{10}}(REM)$ & $0.9329 \pm 0.0003$ & $\log_{{10}}(REM)$ & $0.9096 \pm 0.0003$ \\ 
$\log_{{10}}(MEGNO(LE)/n)$ & $0.93 \pm 0.08$ & $\log_{{10}}(MEGNO(LE)/n)$ & $0.90 \pm 0.09$ \\ 
$\log_{{10}}(\log_{{10}}(LE)/n)$ & $0.924 \pm 0.011$ & $\log_{{10}}(\log_{{10}}(LE)/n)$ & $0.888 \pm 0.015$ \\ 
$\log_{{10}}(FLI^{{WB}}(\hat{{x}}))$ & $0.913 \pm 0.007$ & $\log_{{10}}(FLI^{{WB}}(\hat{{x}}))$ & $0.88 \pm 0.02$ \\ 
$\log_{{10}}(FMA)$ & $0.869 \pm 0.003$ & $\log_{{10}}(FLI(\hat{{x}})/n)$ & $0.843 \pm 0.009$ \\ 
$\log_{{10}}(FLI(\hat{{x}})/n)$ & $0.863 \pm 0.007$ & $\log_{{10}}(FMA)$ & $0.797 \pm 0.005$ \\ 
        \hline
        \hline
        \multicolumn{4}{c}{$(\omega_x, \omega_y) = (0.168, 0.201)$} \\
         \multicolumn{2}{c|}{$\varepsilon = 0.0; \mu = 0.0$} & \multicolumn{2}{c}{$\varepsilon = 0.0; \mu = 0.5$} \\
        \hline
        $\log_{{10}}(GALI^{{(4)}})$ & $0.9896 \pm 0.0004$ & $\log_{{10}}(GALI^{{(4)}})$ & $0.9909 \pm 0.0004$ \\ 
$\log_{{10}}(REM)$ & $0.98682 \pm 0.00009$ & $\log_{{10}}(MEGNO(LE)/n)$ & $0.99 \pm 0.09$ \\ 
$\log_{{10}}(FLI^{{WB}}(\hat{{x}}))$ & $0.986 \pm 0.002$ & $\log_{{10}}(REM)$ & $0.98850 \pm 0.00012$ \\ 
$\log_{{10}}(\log_{{10}}(LE)/n)$ & $0.981 \pm 0.003$ & $\log_{{10}}(FLI^{{WB}}(\hat{{x}}))$ & $0.988 \pm 0.002$ \\ 
$\log_{{10}}(FLI(\hat{{x}})/n)$ & $0.980 \pm 0.016$ & $\log_{{10}}(\log_{{10}}(LE)/n)$ & $0.99 \pm 0.07$ \\ 
$\log_{{10}}(FMA)$ & $0.9319 \pm 0.0010$ & $\log_{{10}}(FLI(\hat{{x}})/n)$ & $0.980 \pm 0.015$ \\ 
$\log_{{10}}(MEGNO(LE)/n)$ & $0.9 \pm 0.2$ & $\log_{{10}}(FMA)$ & $0.9510 \pm 0.0009$ \\ 
        \hline
        \multicolumn{4}{c}{$(\omega_x, \omega_y) = (0.168, 0.201)$} \\
         \multicolumn{2}{c|}{$\varepsilon = 32.0; \mu = 0.0$} & \multicolumn{2}{c}{$\varepsilon = 32.0; \mu = 0.5$} \\
        \hline
        $\log_{{10}}(GALI^{{(4)}})$ & $0.903 \pm 0.003$ & $\log_{{10}}(GALI^{{(4)}})$ & $0.914 \pm 0.003$ \\ 
$\log_{{10}}(REM)$ & $0.8880 \pm 0.0004$ & $\log_{{10}}(REM)$ & $0.8915 \pm 0.0004$ \\ 
$\log_{{10}}(MEGNO(LE)/n)$ & $0.87 \pm 0.09$ & $\log_{{10}}(MEGNO(LE)/n)$ & $0.89 \pm 0.11$ \\ 
$\log_{{10}}(\log_{{10}}(LE)/n)$ & $0.863 \pm 0.016$ & $\log_{{10}}(FMA)$ & $0.881 \pm 0.007$ \\ 
$\log_{{10}}(FLI(\hat{{x}})/n)$ & $0.849 \pm 0.012$ & $\log_{{10}}(\log_{{10}}(LE)/n)$ & $0.88 \pm 0.02$ \\ 
$\log_{{10}}(FLI^{{WB}}(\hat{{x}}))$ & $0.849 \pm 0.007$ & $\log_{{10}}(FLI^{{WB}}(\hat{{x}}))$ & $0.870 \pm 0.012$ \\ 
$\log_{{10}}(FMA)$ & $0.843 \pm 0.004$ & $\log_{{10}}(FLI(\hat{{x}})/n)$ & $0.850 \pm 0.017$ \\ 
        \hline
    \end{tabular}
    \label{table:values}
\end{table*}

Performance estimates have been ranked in decreasing order, separating the various cases considered in our analyses. $GALI^{(4)}$ turns out to be the highest scorer in all cases, followed by $REM$. Then we find $MEGNO$ and $FLI^{{WB}}(\hat{{x}})$, while $FMA$ tends to be the last in this ranking. The error associated with each performance estimate value is provided by the variation of the accuracy whenever the automatic threshold value is varied by $\pm 5\%$. This quantity provides information on the robustness of the accuracy against perturbation of the threshold: A small value indicates a high stability of the numerical values. It is also worth noting that the performance estimates of the best dynamic indicators are correlated with small values of the corresponding error.

\begin{figure*}[htp]
    \centering
    \includegraphics[width=0.85\textwidth]{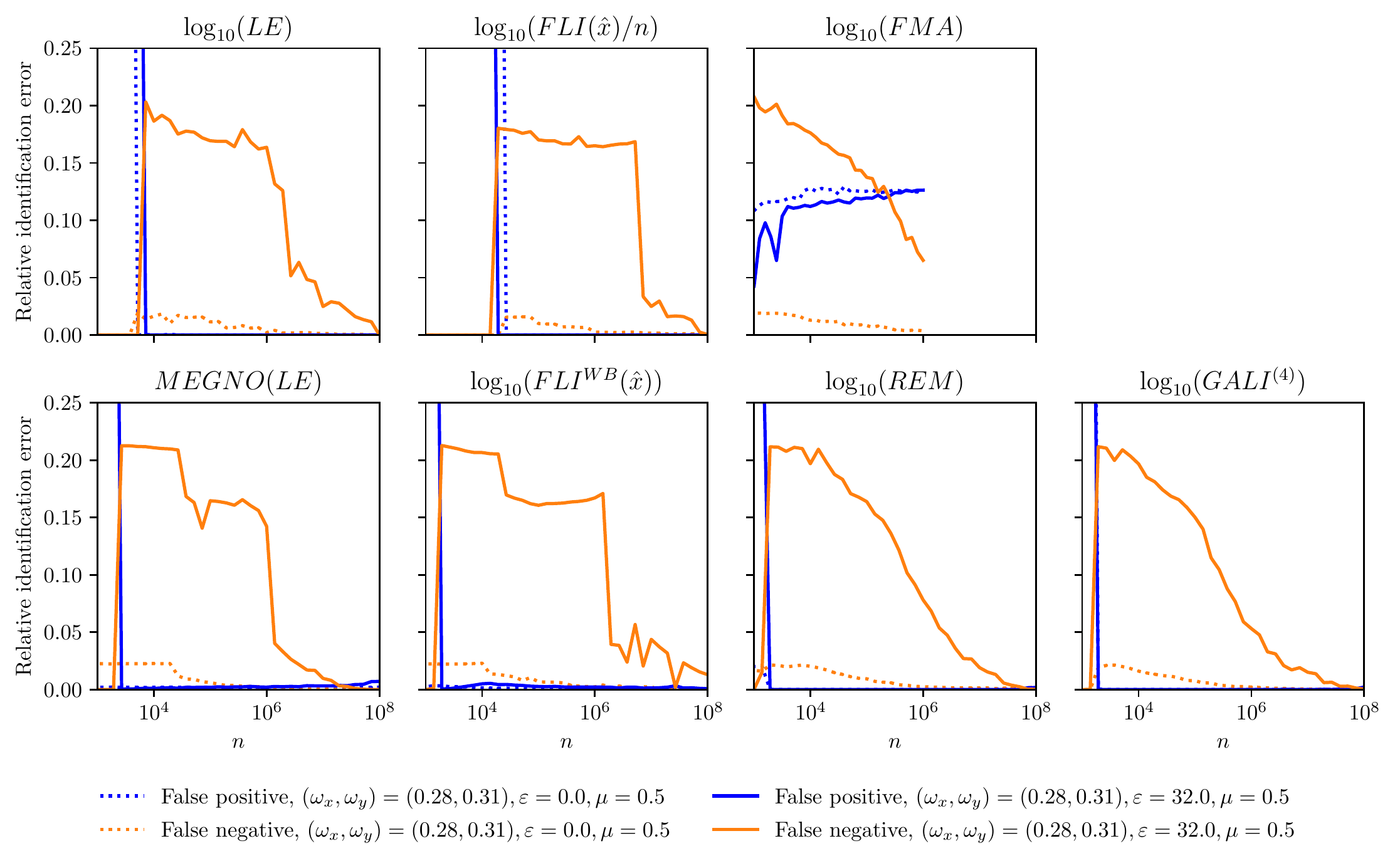}
    \caption{Identification errors for the various indicators as a function of the number of turns for the cases displayed in the first row of Fig.~\ref{fig:performance}.}
    \label{fig:error_comparison}
\end{figure*}

Important insights on the performance of the various indicators can be gained by looking at the relative identification error in terms of false positive, i.e., when a regular orbit is classified as chaotic, and false negative, i.e., when a chaotic orbit is classified as regular. A false negative is almost unavoidable, according to the behavior shown in Fig.~\ref{fig:fli_compare_mean_birk}, unless the indicator is calculated over a very large number of turns, which means accepting a very limited predictive power of the indicator. However, the behavior of the two types of errors reveals interesting features of the various indicators. An overview of the dependence of false positive and false negative errors is shown in Fig.~\ref{fig:error_comparison}, where relative errors are displayed as functions of the turn number for the map configurations considered in the first row of Fig.~\ref{fig:performance}. 

The behavior of the false positive error reveals a fundamental difference between $FMA$ and the other indicators. In fact, $FMA$ shows an error value that is only slightly dependent on the turn number and drops to small values for very large $n$. For the other indicators, for a low number of turns, this type of error is large, and then, around $10^4$ turns, it drops essentially to zero. This feature is related to the fact that, for a low number of turns, the bimodal structure is not yet present. It is also worth noting that, for the case of $FLI(\hat{{x}})$ the Birkhoff averaging introduces a clear improvement by pushing the position of the sudden drop to zero of the false positive error to a lower number of turns.

The false negative error increases sharply at a turn number close to that corresponding to the abrupt decrease in the false positive error. After this turn number, two behaviors are observed: In the first case, the error level is approximately constant until it drops to a low value after $n \approx 10^7$. This value is relatively close to that used to determine the GT, which indicates a limited predictive power of the indicator. In the second case, the error level decreases almost linearly as a function of $n$. This is the key to achieving good performance and is the feature shown by $REM$ and $GALI^{(4)}$. It should be noted that $FMA$ also behaves in this way, i.e., with a linear decrease in the false positive error. However, when the false negative error drops, a jump in the false positive error is observed. This error then shows a decrease that is almost negligible up to $n_\mathrm{max}$. These characteristics, related to the characteristics of the distribution of the $FMA$ values, prevent this indicator from reaching a good performance level. 

As a last comment, these features are always present, but frequency modulation strongly enhances the errors. 
\section{\label{sec:conc} Conclusions}
In this paper, various numerical indicators to identify the chaotic character of orbits of Hamiltonian systems have been presented and discussed in detail. The powerful Birkhoff averages were used to improve the convergence rate of an indicator in the case of regular initial conditions. The goal of our analysis is to evaluate the performance of the indicators in terms of accuracy in the binary classification of an orbit identified by its initial conditions, as regular or chaotic. An important element in this assessment is whether the correct classification can be achieved by using the information over a limited number of turns, i.e., whether an early chaos detection can be effectively performed, which is equivalent to probing the predictive power of dynamic indicators. 

The dynamical system that has been selected as a test bed for performance analyzes is a $4d$ H\'enon-like symplectic map, with or without cubic nonlinearity and with or without frequency modulation. This choice is justified by the relevant applications of this map to understand long-term stability problems in particle accelerators. Several configurations have been considered and, for each case, a ground truth classification has been determined with $n=10^8$ iterations. The various indicators have been used to provide an estimate of the classification performance with respect to ground truth as a function of the number of turns used. The classification is based on the bimodal feature of the indicator value distributions, which points out two clusters associated with regular and chaotic orbits. To define a classification threshold, we use a KDE-based algorithm to determine the position of the distribution minimum between the two modes.

A ranking of the performance of the various indicators has been established, with $GALI^{(4)}$ slightly outperforming the other indicators in all the cases considered, immediately followed by $REM$. Then we find $FLI^{{WB}}$ and $MEGNO(LE)$. Modulation of the linear frequencies significantly reduces the predictive power of each indicator. It should be noted that the identification errors of the various indicators are largely dominated by the wrong labeling of the initial conditions as regular. 
%while $FMA$ also tends to wrongly identify the initial conditions as chaotic. 

The conclusions drawn for the case of the $4d$ H\'enon-like map are generic for a polynomial Hamiltonian system in a neighborhood of elliptic fixed points. Hence, these results can be particularly useful for applications such as nonlinear beam dynamics. The specific choice of an indicator to predict the chaotic character should take into account the performance evaluated in our analysis, as well as the computational efforts needed to compute the various indicators. In this sense, $REM$ could be a very interesting candidate due to its good performance combined with computational efficiency, which is particularly suitable for reducing the CPU time required for the numerical integration of complex physical systems.
\appendix
\section{Computational costs for evaluating the indicators of chaos}\label{app:computing}
Evaluation of a dynamic indicator requires a variable amount of computational cost, which could affect the feasibility and efficiency of specific implementations or favor the usage of specific dynamic indicators. Here, we focus our considerations on the specific case of a discrete map with a known analytic expression for both the tangent and the inverse maps.

For $LE$, $FLI$, and $MEGNO(LE)$, the main computational effort consists of tracking the value of $\Lop_n(\xbf)$, along the orbit of $\xbf$. This implies the additional memory requirement to store a matrix of size $2d\times2d$ and the execution of matrix-matrix and matrix-vector multiplications at each iteration. It should be noted that an important feature of these indicators is that their evaluation at a target iteration number $n$ also provides their value for all lower iteration numbers. This feature frees up additional computational costs for the analysis of the evolution of the dynamic indicator value over time.

$GALI^{(k)}$ requires the evaluation of $\Lop_n(\xbf)$ to calculate the normalized $k$ images of $\boldsymbol{\eta}_j$ with $1\leq j \leq k$. A practical and fast method for computing the norm of external products in Eq.~\eqref{eq2_22} is given in~\cite{Skokos2008}, where it is proven that $GALI^{(k)}$ is equal to the product of singular values $z_j$, of $A$, where $A$ is a $2d\times k$ matrix that reads
\begin{equation}
    A = \left(\begin{array}{ccc}
        \left(\frac{\Lop_n(\xbf) \boldsymbol{\eta}_1}{ \Vert \Lop_n(\xbf) \boldsymbol{\eta}_1\Vert}\right)_1 & \cdots & \left(\frac{\Lop_n(\xbf) \boldsymbol{\eta}_k}{ \Vert \Lop_n(\xbf) \boldsymbol{\eta}_k\Vert}\right)_1 \\
        \vdots &  & \vdots \\
        \left(\frac{\Lop_n(\xbf) \boldsymbol{\eta}_1}{ \Vert \Lop_n(\xbf) \boldsymbol{\eta}_1\Vert}\right)_{2d} & \cdots & \left(\frac{\Lop_n(\xbf) \boldsymbol{\eta}_1}{ \Vert \Lop_n(\xbf) \boldsymbol{\eta}_1\Vert}\right)_{2d}
    \end{array}\right)\,.
\end{equation}

The singular values of $A$ can be obtained by applying the Singular Value Decomposition (SVD) method~\cite{10.5555/1403886}. Note that the evaluation of $GALI^{(k)}$ for a target iteration number $\bar{n}$ also provides the values of $\boldsymbol{\eta}_j$ for all lower values of $n$. However, for each $n \leq \bar{n}$ for which we wish to evaluate $GALI^{(k)}$, a specific SVD calculation is required.

For the reversibility error indicator $BF$, it is possible to use Eq.~\eqref{eq2-12} to evaluate $\boldsymbol{\Xi}_{n}^{BF}(\xbf)$ with the possibility of exploring several realizations of $\vb{\xi}$. This requires the evaluation, for each iteration, of $\Lop_n^{-1}(\xbf)$ or $\Lop_n(\xbf)$, together with the evaluation of the sum with a selected or a set of selected noise realizations. This can lead to higher memory demands when several noise realizations or the time evolution of the indicator needs to be evaluated. Furthermore, its evaluation at a target iteration number $\bar{n}$ does not provide the values for $n \leq \bar{n}$, as each evaluation requires a different summation and noise realization. If the map analyzed is symplectic, the corresponding invariant defined in Eq.~\eqref{eq2_19} can be used, resulting in a computational effort comparable to the evaluation of $LE$.

$REM$, conversely, involves very little computational effort, as it does not require the evaluation of $\Lop_n(\xbf)$, but only the execution of the orbit computation twice. This makes $REM$ very attractive for applications in which no explicit or analytical expression for the tangent map is available. However, the evaluation of $REM$ for a target iteration number $n$ gives no information on its value for lower iteration values, as its evaluation requires separate backtracking each time.

Finally, for $FMA$, if the fundamental frequency is evaluated using FFT-based methods (see, e.g.,~\cite{Bartolini:292773,Bartolini:316949}), considerable effort is required in terms of memory usage, due to the necessity of storing the entire orbit of $M(\vb{x}, n)$, then perform the algorithm. This is not the case if the fundamental frequency is evaluated using the APA method (see, e.g.,~\cite{Bartolini:292773,Bartolini:316949}), as the mean can be progressively evaluated without the need to store the entire orbit history.

Modern parallel computing architectures, such as those offered in General Purpose Graphics Processing Units (GPGPU)~\cite{DBLP:journals/corr/abs-1202-4347}, follow the single-instruction, multiple-data (SIMD) architecture, that is, they execute the same operations over large data allocations, using thread wraps of hundreds of processing cores.

To fully exploit the SIMD architecture, an algorithm must offer options for scaling up parallelization without strong penalties in terms of memory management or branching.

Tracking multiple initial conditions in discrete-time maps is one of the most straightforward processes to implement in a SIMD architecture, as it can be treated as a problem ``embarrassingly parallel''~\cite{Giovannozzi:317866}, and multiple examples of GPGPU applications can be observed, for example, in charged particle tracking in accelerator physics~\cite{pang2014gpu,oeftiger:hb16-mopr025,adelmann2019opal,schwinzerl:ipac21-thpab190,hermes:ipac2022-mopost045,iliakis2022enabling}.

The various indicators of chaos presented here offer, in general, a straightforward conversion to a SIMD approach, since it is immediately possible to perform the tracking and the turn-after-turn dynamic indicator evaluation of several initial conditions. This improvement alone enables mass processing of initial conditions for large values of the turn number $n_\mathrm{max}$, allowing various types of statistical analysis.

However, an exception is given by $FMA$ when evaluated using FFT-based methods, as it requires one to keep track in memory of the orbit of any initial condition and then perform numerical estimates of the fundamental frequencies. Due to this requirement, scaling up the procedure to a large number of turns or a large number of initial conditions may lead to memory limitations. To fully benefit from the SIMD architecture, we evaluated the fundamental frequency via the APA method with Birkhoff weights, which does not require the storage of the entire orbit but only the weighted mean phase advance, which can be progressively evaluated without high memory requirements.

A similar limitation is present in the $BF$ reversibility error, since its direct evaluation, defined in Eq.~\eqref{eq2-12}, requires maintaining track of the entire orbit when there is interest in evaluating different realizations of $\vb{\xi}_n$. In contrast, $REM$ offers a straightforward GPGPU approach, since it only requires explicit forward and backward tracking, without the need to evaluate the tangent map. We recall that $REM$ evaluates only the first invariant from a single noise realization, obtained by exploiting the numerical roundoff.
\section{Time dependence of dynamic indicators\label{app:timedep}}

\begin{figure*}[ht]
    \centering
    \includegraphics[width=0.9\textwidth]{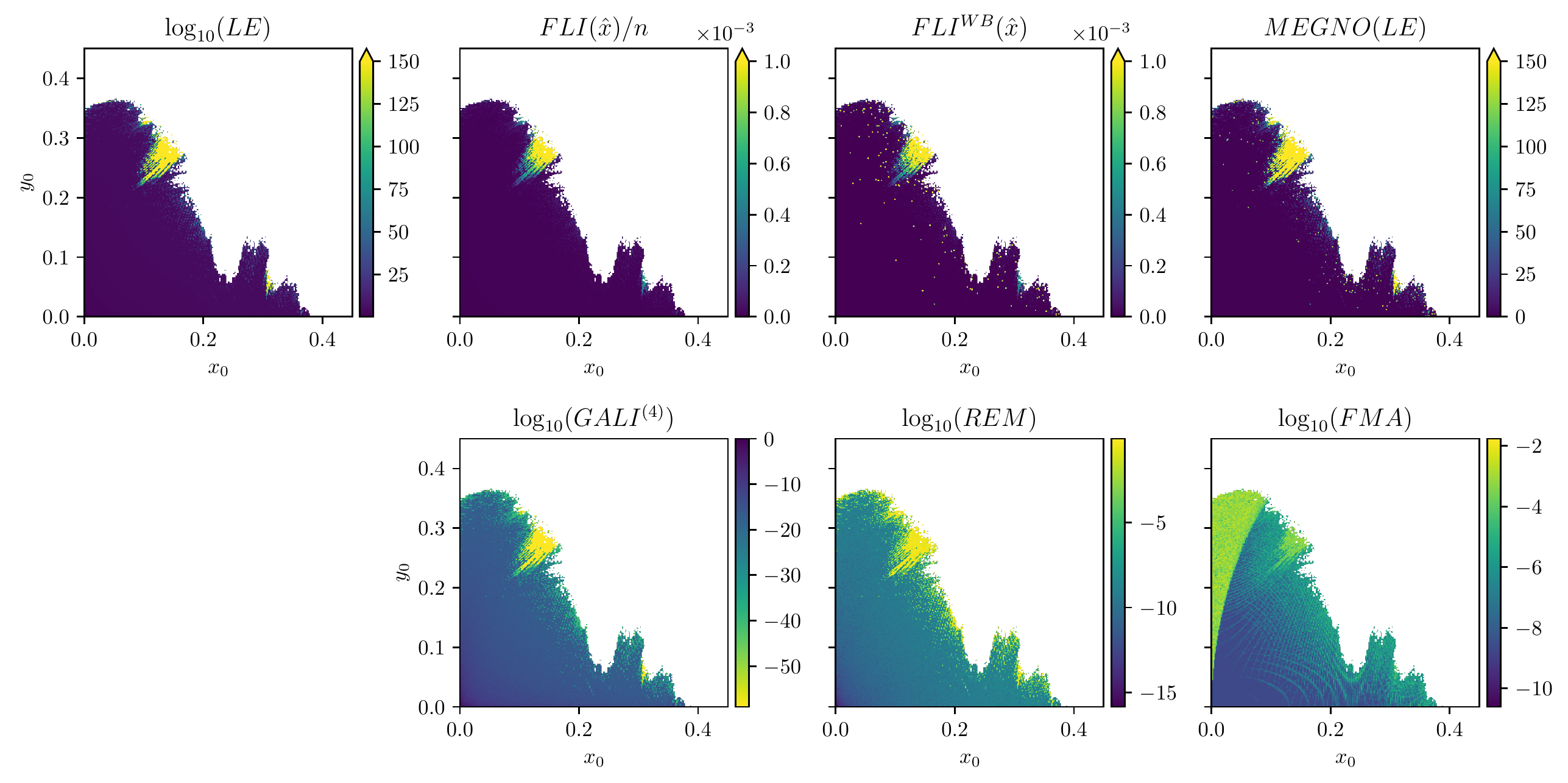}
    \caption{Color maps of the various dynamic indicators for a modulated $4d$ Hénon map evaluated at $n=10^5$. It can be seen how the indicators globally highlight the same structures in phase space, with the exception of $FMA$, which also shows structures related to resonances. Note that an arrow at the top of the color bar means that pixels of the top color correspond to a value equal to or greater than the top value. White pixels correspond to initial conditions whose distance from the origin has exceeded a predefined radius ($r_c=10^2$) during the tracking, before reaching the target iteration number $n_\text{max}=10^8$. (Simulation parameters: $(\omega_{x0},\omega_{y0})= (0.28,\ 0.31),\ \varepsilon=32.0,\ \mu=0.5$).}
    \label{fig:generic_example}
\end{figure*}

\begin{figure*}[ht]
    \includegraphics[width=0.9\textwidth]{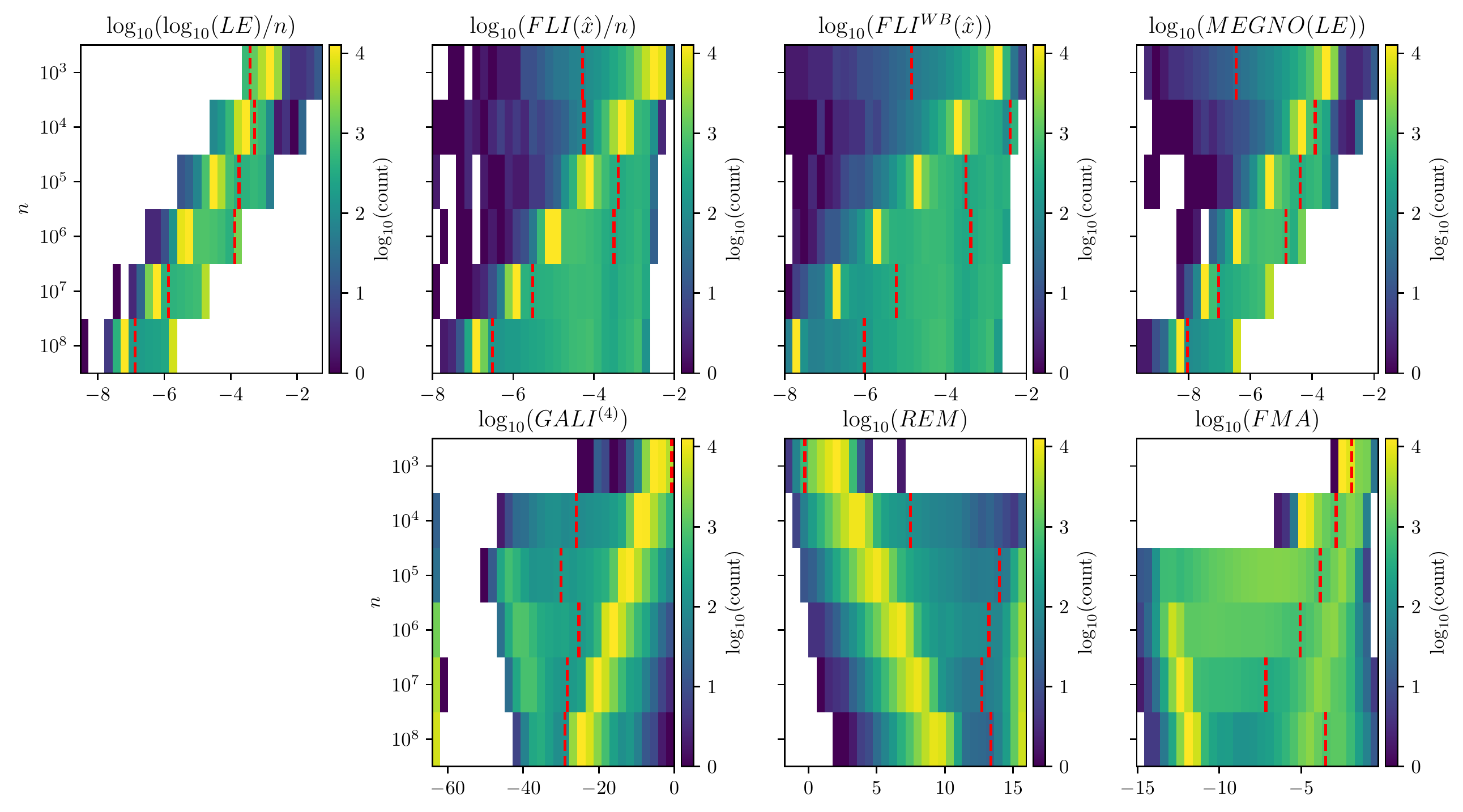}
    \caption{Distribution of values of the various dynamic indicators as a function of time for a modulated $4d$ Hénon map. For low values of the  iterations $n$, the distribution is in general represented by a uni-modal function. For higher values of $n$, we can see the formation of either two separate clusters, making the distribution bi-modal, or an individual cluster with a significant tail. $\log_{10}(FMA)$ constitutes an exception, as it evolves forming a tri-modal distribution (also shown in detail in Fig.~\ref{fig:thresholds}, bottom). The red dashed lines represent threshold values, defined by our algorithm, representing our criterion to distinguish regular and chaotic orbits. (simulation parameters: $(\omega_{x0},\omega_{y0})= (0.28,\ 0.31),\ \varepsilon=32.0,\ \mu=0.5$).}
    \label{fig:generic_example_2}
\end{figure*}

When considering a large amount of initial conditions to determine the properties of the corresponding orbits by means of dynamic indicators, it is possible to obtain an accurate picture of the phase-space structures, such as regions characterized by regular dynamics and regions where frequency modulation and nonlinearities induce chaotic behavior. In Fig.~\ref{fig:generic_example}, the seven chaos indicators computed for $n=10^5$ are presented for a set of initial conditions that turned out to be stable up to $n_\text{max}=10^8$. All indicators highlight a region of regular motion close to the origin and chaotic structures at higher amplitudes. Generally speaking, the various dynamic indicators reconstruct very similar shapes for the regular and chaotic regions of the phase space, with the exception of $FMA$. Indeed, this indicator provides a lot of structure even inside the region that is classified as regular by the other indicators, and in which the values of the other indicators are to a high degree of accuracy constant. We inspect the distribution of values of the various dynamic indicators, computed at a large number of turns. It is possible to observe the formation of bimodal or, as we shall see for the case of $FMA$, three-modal distributions. In Fig.~\ref{fig:generic_example_2}, the time evolution of the distribution of the indicator value is shown. The red lines represent the threshold that we use to distinguish between regular and chaotic orbits, whose definition was given in Section~\ref{subsec:classification}.

The indicators $\log_{10}(LE)/n$, $FLI(\hat{x})/n$, $FLI^{WB}(\hat{x})$, and $MEGNO(LE)/n$ have a comparable behavior and globally tend to cluster regular orbits into an ensemble peaked at near-zero values, whereas chaotic orbits are part of another cluster featuring a large spread of values, which correspond to indicator values that are orders of magnitude higher. To achieve a valid overview of the value distribution, especially its tendency to create a bimodal distribution~\cite{PhysRevE.60.2761,VALLEJO200326}, we will consider the logarithm of these three indicators, i.e.~$\log_{10}(\log_{10}(LE)/n)$, $\log_{10}(FLI(\hat{x})/n)$, $\log_{10}(FLI^{WB}(\hat{x}))$, and $\log_{10}(MEGNO(LE)/n)$.

By direct inspection of the color maps in Fig.~\ref{fig:generic_example_2}, it can be seen how these four indicators generate bimodal distributions, with the peak corresponding to regular orbits featuring a clear trend towards zero, and this trend appears to be faster for $\log_{10}(FLI^{WB}(\hat{x}))$, and $\log_{10}(MEGNO(LE)/n)$, due to the applied filters. $\log_{10}(FLI(\hat{x})/n)$ and $\log_{10}(FLI^{WB}(\hat{x}))$ feature an increasing spread of values corresponding to chaotic orbits, a clear trend of the distribution of regular orbits toward zero. A similar trend is also observed in $\log_{10}(\log_{10}(LE)/n)$ and $\log_{10}(MEGNO(LE)/n)$, however, the current numerical implementation of $LE$ suffers from numerical saturation for chaotic orbits that exhibit exponential growth in the values of the tangent map. This results in a limitation for the spread of values that can be observed for chaotic orbits at high numbers of turns, but, ultimately, the distinction between clusters remains. 

$GALI^{(4)}$ takes values in the interval $[0,1]$, corresponding to the range of values of the volume of the $4d$ parallelotope, constructed by normalized orthonormal displacements. The unit value is associated with the initial orthonormal displacement, whereas zero implies an exact chaos-induced alignment of at least two displacement vectors along the direction of the maximum Lyapunov exponent. Inspecting the indicator distribution in logarithmic scale, i.e.~$\log_{10}(GALI^{(4)})$, highlights a bimodal distribution, where the peak corresponding to the ensemble of regular orbits moves towards small values of the indicator, following a power law distribution. Moreover, an ensemble of chaotic orbits creates a tail distribution of values lower than the regular ensemble, thus creating a second, smaller-amplitude peak in the indicator distribution. The presence of the logarithm when evaluating the distribution of $GALI^{(4)}$ generates a numerical artifact. Indeed, certain chaotic orbits feature a $4d$ volume, computed using the SVD method, that reaches values below numerical precision, which are consequently registered as zero. We assign to these initial conditions a value of $10^{-64}$, which represents a product of 4 singular values $z_j=\epsilon\sim10^{-16}$ with extended precision. The cluster of these special initial conditions generates yet another peak in the indicator distribution that is, nevertheless, irrelevant in future considerations about the classification of orbits. 

The dynamic indicator $REM$ is also considered on a logarithmic scale to better appreciate its behavior. The measured Euclidean distance for the case of regular orbits ranges from a few orders of magnitude higher than the numerical precision $\epsilon\sim 10^{-16}$ for small values of $n$. These indicator values increase with $n$ following a power law (typically, the peak reaches $~10^{5}$ for $n=10^5$) due to the accumulation of the numerical error. Instead, for chaotic orbits, we observe exponential growth that saturates to an almost constant value. This occurs since chaotic orbits belong to an invariant bounded set of diameter $D$ so that the saturation value is about $\epsilon^{-1} D$. %The origin of such an upper bound in the value of $REM$ is found in the fact that each orbit considered for the evaluation of the indicator is limited within a finite-radius sphere.
Similarly to $GALI^{(4)}$, we inspect the indicator in logarithmic scale, i.e.~$\log_{10}(REM)$. 

$FMA$ is based on the evaluation of the Euclidean distance in the frequency space of the fundamental frequencies computed over different time intervals. If we inspect its distribution on logarithmic scale, we observe how the indicator converges to a three-mode distribution. This configuration consists of an ensemble of initial conditions rapidly converging to values close to numerical precision, an ensemble of initial conditions maintaining values above $10^{-5}$, and a well-populated ensemble of initial conditions that connect these two ensembles (this distribution is also shown in Fig.~\ref{fig:thresholds}, bottom). Inspecting the logarithm of the indicator, i.e.~$\log_{10}(FMA)$, allows to inspect the full spread of values achieved by the various orbits.

\bibliographystyle{unsrt}
\bibliography{indicators}

\end{document}